\immediate\write18{makeindex -s nomencl.ist -o "\jobname.nls" "\jobname.nlo"}

\documentclass[journal,onecolumn,draftclsnofoot,]{elsarticle}

\usepackage{lineno,hyperref}
\modulolinenumbers[254]
\usepackage[ngerman,UKenglish]{babel}
\usepackage[top=55pt,bottom=60pt,left=60pt,right=55pt]{geometry}

\usepackage{multicol}
\usepackage{multirow}
\usepackage{arydshln}
\usepackage{lipsum}
\usepackage{algorithm}
\usepackage{algorithmic}

\usepackage{setspace}
\usepackage{amsmath,bm}
\usepackage{fixmath}
\usepackage{commath}
\usepackage{lipsum}
\usepackage{textcomp}

\usepackage[bitstream-charter]{mathdesign}
\usepackage[T1]{fontenc}
\usepackage[utf8]{inputenc} %unicode support

\usepackage{framed}
\usepackage{nomencl} % Nomenclature package
\usepackage{tikz}
\usepackage{colortbl}

\usepackage{ifthen} % For conditional formatting

\renewcommand{\nomgroup}[1]{%
   \ifthenelse{\equal{#1}{S}}{\item[\textbf{Symbols}]}{% Symbols group
   \ifthenelse{\equal{#1}{A}}{\item[\textbf{Abbreviations}]}{% Abbreviations group
   }}}

\makenomenclature % Activate nomenclature    

\biboptions{sort&compress} %To make [1,2,3] be [1-3]

\usepackage{calrsfs}
\DeclareMathAlphabet{\pazocal}{OMS}{zplm}{m}{n}
\let\mathcal\undefined
\newcommand{\mathcal}[1]{\pazocal{#1}}

\usepackage{pifont}

\usepackage{comment}
\usepackage{mathtools}
\usepackage[font=small,skip=6pt]{caption}

\usepackage{xcolor,colortbl}
\usepackage{float}
\usepackage{xcolor}
\usepackage{booktabs}
\definecolor{lightgray}{gray}{0.95}
\definecolor{darkgreen}{rgb}{0.0, 0.2, 0.13}
 \definecolor{darkorchid}{rgb}{0.6, 0.2, 0.7}
\usepackage{dblfloatfix} 
\usepackage{placeins}
\usepackage[font=small,skip=6pt]{caption}

\usepackage[symbol]{footmisc}
\def\correspondingauthor{\footnote{$^*$Corresponding author.}}

\usepackage[utf8]{inputenc}
\usepackage{framed}
\usepackage{nomencl} 
\makenomenclature
\renewcommand*\nompreamble{\begin{multicols}{1}}
\renewcommand*\nompostamble{\end{multicols}}

\journal{Journal of \LaTeX\ Templates}

\journal{XXXX}

\usepackage{cleveref} 

\usepackage{hyperref}
\hypersetup{
    colorlinks=true,
   linkcolor=blue,
    citecolor=black,
    pdftitle={PhyML for fatigue lifetime prediction},
}

\let\oldhref\href
\renewcommand{\href}[2]{\oldhref{#1}{\hbox{#2}}}

\begin{document}

\begin{frontmatter}

\title{Physics-based machine learning for fatigue lifetime prediction\\ under non-uniform loading scenarios}

\author{Abedulgader Baktheer}
%\ead{abedulgader.baktheer@ibnm.uni-hannovr.de}
\author{Fadi Aldakheel\correspondingauthor{}$^*$}
\ead{fadi.aldakheel@ibnm.uni-hannovr.de}
\ead[url]{https://www.ibnm.uni-hannover.de}

\address{Institute of Mechanics and Computational Mechanics, Leibniz Universität Hannover, 30167 Hannover, Germany}

\begin{abstract}

Accurate lifetime prediction of structures and structural components subjected to cyclic loading is vital, especially in scenarios involving non-uniform loading histories where load sequencing critically influences structural durability. Addressing this complexity requires advanced modeling approaches capable of capturing the intricate relationship between loading sequences and fatigue lifetime. Traditional high-cycle fatigue simulations are computationally prohibitive, necessitating more efficient methods.
This work highlights the potential of physics-based machine learning ($\phi$ML) to predict the fatigue lifetime of materials under various loading conditions. Specifically, a feedforward neural network is designed to embed physical constraints from experimental evidence, including initial and boundary conditions, directly into its architecture to enhance prediction accuracy. It is trained using numerical simulations generated by a physically based anisotropic continuum damage fatigue model.
The model is calibrated and validated against experimental fatigue data of concrete cylinder specimens tested in uniaxial compression. The simulations used for training quantify the effects of load sequences considering scenarios under two different loading ranges.
The proposed approach demonstrates superior accuracy compared to purely data-driven neural networks, particularly in situations with limited training data, achieving realistic predictions of damage accumulation. To this end, a general algorithm is developed and successfully applied to predict fatigue lifetimes under complex loading scenarios with multiple loading ranges. In this approach, the $\phi$ML model serves as a surrogate to capture damage evolution across load transitions. The $\phi$ML based algorithm is subsequently employed to investigate the influence of multiple loading transitions on accumulated fatigue life, and its predictions align with trends observed in recent experimental studies.
The presented contribution demonstrates physics-based machine learning as a promising technique for efficient and reliable fatigue life prediction in engineering structures, with possible integration into digital twin models for real-time assessment.

\end{abstract}

%Loading sequence effect, Variable amplitudes

\begin{keyword}
$\phi$ML, Anisotropic fatigue model, Experimental studies, Lifetime prediction, Digital twin
\end{keyword}

\end{frontmatter}

%\linenumbers

\section{Introduction}
\label{sec:introduction}
The lifespan of structures and structural components is typically governed by the progression of damage effects throughout their service life. Fatigue damage accumulation is a key factor in the design and assessment of structures subjected to cyclic loading, and it has gained significant attention in recent years~\cite{li2025arresting, Baktheer_2024_stress, chen2024effect, leuders2013mechanical, Ali_2022, van2022improving, MIARKA_2022, Seles_2021}. This is driven by the growing demand for highly efficient, lightweight engineering infrastructures in energy and transportation sectors, such as wind power plants and bridges, while simultaneously addressing the challenge of mitigating early fatigue failures in aging structures.
Current design codes and standards provide simplified methodologies to estimate accumulated fatigue damage such as the Palmgren-Miner (P-M) rule \cite{palmgren1924lebensdauer, miner1945cumulative}, relying primarily on constant-amplitude fatigue test data \cite{fib_model_code_fib_2010, EN19922, EN19932, EN19931-9} (e.g., S-N curves). However, real-world loading conditions often deviate substantially from constant-amplitude loading, highlighting a critical need to deepen our understanding of fatigue behavior under more realistic and complex loading scenarios.
Exclusively addressing this challenge through experimental investigations is impractical due to the vast range of potential load combinations and conditions. Consequently, advanced modeling approaches capable of capturing the intricate interactions between loading sequences and fatigue lifetime are essential for overcoming these limitations.

\begin{table*}[!t]   
%\vspace{-2mm}
\begin{framed}
\small
\nomenclature[S]{$\eta_\mathrm{i}$}{consumed fatigue lifetime in the range i}
\nomenclature[S]{$N_\mathrm{i}$}{number of loading cycles of the range i}
\nomenclature[S]{$N_\mathrm{i}^\mathrm{f}$}{the number of cycles to failure under a constant \newline amplitude of the range i}
\nomenclature[S]{$\phi$}{free energy potential }
\nomenclature[S]{$\boldsymbol{\varepsilon}$}{strain tensor}
\nomenclature[S]{$\boldsymbol{\omega}$}{damage tensor}
\nomenclature[S]{$\lambda$}{first Lam$\Acute{\mathrm{e}}$ constant}
\nomenclature[S]{$\mu$}{second Lam$\Acute{\mathrm{e}}$ constant}
\nomenclature[S]{$\mathrm{tr}$}{trace of the tensor}
\nomenclature[S]{$g$}{material parameter of the applied fatigue model}
\nomenclature[S]{$\alpha$}{material parameter of the applied fatigue model}
\nomenclature[S]{$\beta$}{material parameter of the applied fatigue model}
\nomenclature[S]{$\boldsymbol{I}$}{identity tensor}
\nomenclature[S]{$\boldsymbol{Y}$}{energy release rate tensor}
\nomenclature[S]{$f(\boldsymbol{\varepsilon},\boldsymbol{\omega} )$}{yield function of the applied fatigue model}
\nomenclature[S]{$k(\boldsymbol{\omega})$}{damage threshold}
\nomenclature[S]{$C_0$}{material parameter of the applied fatigue model}
\nomenclature[S]{$C_1$}{material parameter of the applied fatigue model}
\nomenclature[S]{$(.)^{+}$}{positive part of a tensor}
\nomenclature[S]{$I_2$}{second invariant of the energy release rate tensor}
\nomenclature[S]{$\tilde{\varepsilon}$}{equivalent strain}
\nomenclature[S]{$f$}{material parameter of the applied fatigue model}
\nomenclature[S]{$n$}{material parameter of the applied fatigue model}
\nomenclature[S]{$K$}{material parameter of the applied fatigue model}
\nomenclature[S]{$S_\mathrm{i}^\mathrm{max}$}{the upper level of the loading range i}
\nomenclature[S]{$S$}{loading ratio}
\nomenclature[S]{$S_\mathrm{i}^\mathrm{min}$}{the lower level of the loading range i}
\nomenclature[S]{$f_{\mathrm{c}}$}{compressive strength of the concrete}
\nomenclature[S]{$\omega$}{damage parameter}
\nomenclature[S]{$\Delta S^\mathrm{max}$}{upper jump between two loading ranges}
\nomenclature[S]{$\Delta S^\mathrm{min}$}{lower jump between two loading ranges}
\nomenclature[S]{$\bar{S}_\mathrm{i}$}{mean loading level within two loading ranges}
\nomenclature[S]{$S_\mathrm{i}^\mathrm{m}$}{mean value of individual loading range}

\nomenclature[S]{$\eta^\mathrm{out}_\mathrm{i}$}{predicted sum of consumed fatigue lifetime at \newline each load jump}
\nomenclature[S]{$\eta^\mathrm{new}_\mathrm{i}$}{corrected  consumed fatigue lifetime before each \newline new load jump}
\nomenclature[S]{$\sum \eta $}{sum of consumed and remaining fatigue lifetime \newline (accumulated fatigue lifetime)}
\nomenclature[S]{$\eta^\mathrm{rem}$}{remaining fatigue lifetime}
\nomenclature[S]{$\Delta \eta_\mathrm{i}$}{correction of consumed fatigue lifetime \newline in comparison to P-M rule}
%\nomenclature[S]{$\sum \eta_\mathrm{i}$}{corrected  consumed fatigue lifetime }
\nomenclature[S]{$\eta^\mathrm{cons}$}{consumed fatigue lifetime before the loading jump}

\nomenclature[S]{$\mathcal{L}_{\text{cust}}$}{custom loss function of the $\phi$ML model}
\nomenclature[S]{$\mathcal{L}_{\text{data}}$}{data based loss function}
\nomenclature[S]{$\mathcal{L}_{\text{const}}$}{physical constraints based loss function}
\nomenclature[S]{$\mathcal{L}_{\text{bound}}$}{boundary and initial \newline conditions based loss function}
\nomenclature[S]{$\mathcal{L}_{\text{spars}}$}{sparse domain regularization based loss function}
\nomenclature[S]{$\mathcal{L}_{\text{total}}$}{total loss of the $\phi$ML model}

% \nomenclature[S]{$w_{\text{data}}$}{weight parameter of the data loss}
\nomenclature[S]{$w_{\text{const}}$}{weight parameter of the physical constraints loss}
\nomenclature[S]{$w_{\text{bound}}$}{weight parameter of the boundary and initial conditions loss}
\nomenclature[S]{$w_{\text{spars}}$}{weight parameter of the sparse domain regularization loss}

\nomenclature[S]{$y_{\text{pred}}$}{predicted output value by the machine learning model}
\nomenclature[S]{$y_{\text{target}}$}{ground truth value}
\nomenclature[S]{$\epsilon$}{loss tolerance threshold}
\nomenclature[S]{$\zeta$}{constant of the ELU activation function}

%%%%%%%%%%%%%%%%%%Abbreviations%%%%%%%%%%%%%%
%\nomenclature[A]{NN}{Neural Network}
\nomenclature[A]{FFNN}{Feed-Forward Neural Network}
%\nomenclature[A]{$\phi$NN}{Physics-based Neural Network}
\nomenclature[A]{$\phi$ML}{Physics-based Machine Learning}
\nomenclature[A]{H-L}{High-Low loading scenario}
\nomenclature[A]{L-H}{Low-High loading scenario}
\nomenclature[A]{P-M}{Palmgren-Miner rule}
\nomenclature[A]{ML}{Machine Learning}
\nomenclature[A]{LS}{Loading Scenario}
\nomenclature[A]{ELU}{Exponential Linear Unit}
\nomenclature[A]{MSE}{Mean Squared Error}

\printnomenclature

\end{framed}
%\vspace{-4mm}
\end{table*}

Several fatigue modeling approaches have been developed over the last few years for various materials, based on damage mechanics, phase field modeling, and plasticity frameworks~\cite{marigo1985modelling, lemaitre2005engineering, caggiano_2020, polym_6020, ULLOA_2021, ALESSI_2018, Carrara_2020, Ghafoori_2020, schreiber_2021, KHALIL_2022, SCHRODER_2022}.
Multiscale models have been refined to capture the unique microcrack-induced stress redistribution characteristics under cyclic loading~\cite{kobler2021computational, lee2021review, kirane2015microplane, BAKTHEER_2021_MS1, Baktheer_2024_FFEMS, hessman2023micromechanical, Ueda2019, LE_2023, shojai2024micro}. However, while these approaches provide detailed insights, they remain {\it computationally intensive}, particularly for high-cycle fatigue simulations.
Given the high computational cost of performing fatigue simulations for large components and structures, there is a pressing need for innovative approaches and efficient surrogate models. These advancements will be essential to translate the expected growth in fatigue modeling into practical and efficient lifetime-based simulations for real-world applications.

\subsection{Machine learning methods}
Classical numerical methods for fatigue analysis are struggling to keep up with the growing complexity of modern engineering systems, both mathematically and computationally. In this regard, machine learning (ML) has emerged as a promising tool for advancing fatigue lifetime prediction. By leveraging its ability to recognize and learn intricate patterns from data, ML offers new opportunities for improving the efficiency and accuracy of fatigue modeling.

Recent efforts have focused on developing ML models to predict fatigue lifetime based on experimental datasets, addressing the challenges of nonlinear behaviors and complex loading conditions, e.g.,~\cite{zhan2021machine, Son_2022, Zhang_2022, zhan2021data, RIYAR_2023, hamada2025advancing}. Data-driven approaches, such as back-propagation neural networks, have demonstrated potential in predicting fatigue crack growth rates and fatigue life, effectively leveraging experimental datasets while reducing data acquisition costs through active learning~\cite{ZHANG_2022_ML_fatigue_1}. Advanced deep learning mechanisms, like self-attention, have been incorporated to accurately capture the effects of nonlinear loading histories and varying temperatures on multiaxial fatigue life~\cite{YANG_2022_ML_fatigue_2}.

While these black-box ML methods are often more efficient than traditional numerical methods, their accuracy depends largely on the amount of available training data~\cite{app15052589}. This can make them unreliable, unrobust, and uninterpretable, i.e., their accuracy cannot be guaranteed, they are sensitive to changes in the data they are trained on, and it is difficult to understand how they make predictions. Thus, they do not generalize well when evaluated outside the range of their training data~\cite{montans2023}. All of these aspects are critical in fatigue problems, where, on the one hand, only little experimental data are available or are expensive to generate, and, on the other hand, high demands on the accuracy and reliability of models and simulation results are required due to safety and certification regulations.

\subsection{Physics-based machine learning}

To tackle the aforementioned challenges, only in the last 4 to 5 years researchers {from mathematics,  computer, and engineering sciences} have begun to {systematically combine the well-known principles of physical modeling and established numerical simulation methods with ML}, i.e., started developing physics-based machine learning ($\phi$ML) methods.
In the literature, such approaches are also referred to as scientific, structure-preserving or hybrid machine learning \cite{karniadakis2021piml, ASAD_2023, michopoulos2024scientific, Cueto_2023, ELFALLAKIIDRISSI_2024}. Most prominently, {physics-informed} neural networks (PINNs) incorporate residual terms of differential equations into the loss function of the NN \cite{di2024physics, manav2024phase, THIERCELIN_2024, LEON_2025, Athanasiou_2025}.

For fatigue failure, Salvati et al.~\cite{SALVATI_2022_ML_fatigue_3} developed a PINN framework that combines machine learning with phenomenological laws to improve finite-life fatigue prediction in materials with defects. The approach incorporates a linear elastic fracture mechanics based semi-empirical model and is validated using experimental data from additively manufactured materials. Other hybrid approaches, such as integrating Long Short-Term Memory cells with physical laws, further enhance prediction accuracy while maintaining physical consistency~\cite{HALAMKA_2023_ML_fatigue_6, HE_2023__ML_fatigue_8}. Wang et al.~\cite{WANG_2023_ML_fatigue_7} proposed two physics-guided machine learning frameworks for predicting the fatigue life of additive manufacturing materials. The first is a hybrid ML model that incorporates morphological details of critical defects to enhance fatigue life prediction. The second, a serial physics-guided ML framework, demonstrates that artificial neural networks achieve higher prediction accuracy than support vector regression. For elevated temperature fatigue, PINN enhances accuracy by enforcing material constraints, ensuring experimental consistency and efficiency~\cite{JIANG_2024_ML_fatigue_5}.

Despite these advancements, challenges such as computational inefficiency, sensitivity to hyperparameters, and difficulties in ensuring physical consistency highlight the need for a {\it stronger integration} of physical constraints into machine learning frameworks. Addressing this critical gap is essential for developing realistic and practical approaches to fatigue lifetime prediction under various loading conditions. This motivates the development of physics-based machine learning ($\phi$ML), to overcome these challenges. So far, $\phi$ML models that aim to incorporate thermodynamic consistency have been proposed for hyperelasticity~\cite{gonzalez2020,klein2022}, as well as inelastic behaviors such as visco-elasticity \cite{tac2023a,rosenkranz2024a} or elasto-plasticity~\cite{vlassis2021,masi2023,fuhg2023}, brittle and ductile fracture mechanics \cite{Aldakheel_2025}. However, fully \emph{physics-augmented} inelastic formulations are yet to be extended to fatigue lifetime prediction.

\subsection{Experimental background} \label{sec:experimental_background}
The characterization of fatigue under non-uniform loading conditions presents a significant challenge, leading to a substantial gap in research that remains largely unexplored across various materials. Understanding how loading sequences influence fatigue behavior is crucial, yet only a limited number of experimental studies have addressed this issue~\cite{Agerskov}. For instance, in the case of concrete, research on this aspect remains limited~\cite{Tepfers_1977}, with most studies focusing on elementary scenarios such as simple two-stage high-low (H-L) or low-high (L-H) cyclic loading. While these approaches provide valuable insights, they fall short of addressing the more complex loading conditions encountered in practice.

Experimental investigations reveal that loading sequences significantly influence fatigue lifetime. In compression, the H-L sequence generally reduces fatigue lifetime compared to the L-H sequence~\cite{holmen1982fatigue, petkovic1990fatigue, Baktheer_2021_3,  hoff1984testing}. This contradicts the predictions of the widely used Palmgren-Miner (P-M) rule, which assumes linear damage accumulation based solely on constant-amplitude fatigue test data. In tension, however, limited experimental evidence indicates a reverse trend, where the L-H sequence reduces fatigue lifetime more than the H-L sequence~\cite{hilsdorf1966fatigue}.
Similar observations have been recorded in the experiments monitored with digital image correlation in~\citep{BAKTHEER_2021_4}.
Further studies on the bond fatigue between concrete and steel reinforcement confirm these findings~\cite{ baktheer2018modeling}.
The simplicity of the P-M rule, while convenient for design codes, fails to account for these sequence-dependent effects. This underscores the need for more advanced investigations through modeling approaches to capture the complexities of real-world loading scenarios.

\subsection{Objective of the Study}
The current paper aims to explore the potential of $\phi$ML (by integrating physical knowledge and constraints into artificial neural networks) to enable efficient prediction of the fatigue lifetime of high-strength concrete, explicitly considering the effects of loading sequence in non-uniform loading scenarios. 
An efficient surrogate model is developed by training a physics-based machine learning framework using fatigue simulations generated from a representative macro-scale anisotropic continuum damage model.
The resulting surrogate model offers two key applications:
\begin{itemize}
    \item As an advanced tool for estimating cumulative fatigue damage under variable amplitude loading scenarios. It can be directly integrated into digital twin models for real-time assessment in engineering structures, eliminating the need for computationally expensive nonlinear fatigue simulations.
    \item Integration with finite element simulations to efficiently model fatigue behavior over the lifetime of a structure. By accounting for the number of cycles and corresponding stress amplitudes at each material point, it enables precise evaluation of stiffness degradation after individual or multiple load cycles.
\end{itemize}
In this work, we focus on the first approach, developing a surrogate model for fatigue damage estimation within digital twin frameworks, offering an efficient and practical alternative to nonlinear fatigue simulations.

%\textcolor{red}{Structure of the paper:}
The structure of this paper is organized as follows: Section~\ref{sec:fatigue_model} provides a summary of the anisotropic continuum damage fatigue model used to generate the training data for the physics-based machine learning framework, along with details on its calibration and validation.
Section~\ref{sec:neural_networks} describes the development of the $\phi$ML model, including the data generation process, the chosen architecture, the integration of physical constraints, the customized loss function, and the training procedure.
In Section~\ref{sec:results_discussion}, the results obtained using the $\phi$ML model are presented and discussed. Comparative analyses are performed between purely data-driven models and those incorporating physical constraints, considering both large and small datasets.
Section~\ref{sec:generalization_validation} extends the developed $\phi$ML framework to accommodate arbitrary loading scenarios through a proposed algorithm. This section also includes experimental validation studies of fatigue lifetime predictions under more complex loading conditions.

\section{Anisotropic continuum damage model}
\label{sec:fatigue_model}
The fatigue damage modeling approach developed by \cite{HALM1998439, DRAGON2000331, alliche2004damage} is selected for this study due to its computational efficiency. Its incremental formulation ensures compatibility with a wide range of fatigue loading histories, making it adaptable to complex scenarios. The model employs a stress-driven formulation that simplifies the problem to a uniaxial stress state, allowing efficient simulation of compressive loading cycles on cylindrical solids under the assumption of uniform stress distribution. This computational efficiency is critical for conducting the large number of fatigue simulations needed to train the machine learning model. Furthermore, the model effectively captures the nonlinear progression of fatigue damage and incorporates the influence of loading sequence~\cite{BAKTHEER2019}. This is achieved by replacing the conventional yield limit with an irreversibility condition, ensuring accurate representation of loading and unloading behavior.
It is important to highlight that the model used in this study primarily serves as a demonstrative example to showcase the potential of physics-based machine learning. Specifically, it illustrates how advanced fatigue damage modeling knowledge can be effectively transferred into ML frameworks for efficient lifetime prediction. While this model was chosen for its simplicity and adaptability, other advanced fatigue modeling approaches \cite{BAKTHEER_2021_MS1, Seles_2021,ULLOA_2021,Carrara_2020}, could equally serve as the foundation for similar applications.

\subsection{Model formulation}
The model is based on the anisotropic damage framework, defined by the following free energy potential
\begin{align}\label{eq:Helmholtz_free_energy}
\phi(\boldsymbol{\varepsilon}, \boldsymbol{\omega})  =   \frac{1}{2} \lambda\;[\mathrm{tr} (\boldsymbol{\varepsilon})]^2  + \mu \; \mathrm{tr}(\boldsymbol{\varepsilon} \cdot \boldsymbol{\varepsilon})  + g \; \mathrm{tr}(\boldsymbol{\varepsilon} \cdot \boldsymbol{\omega}) 
+\alpha \; \mathrm{tr}(\boldsymbol{\varepsilon}) \;\mathrm{tr}(\boldsymbol{\varepsilon} \cdot \boldsymbol{\omega}) + 2\;\beta \; \mathrm{tr}(\boldsymbol{\varepsilon} \cdot \boldsymbol{\varepsilon} \cdot \boldsymbol{\omega}) \; ,
\end{align}
where $\lambda, \mu$ are the Lam$\Acute{\mathrm{e}}$ constants and $\alpha, \beta, g$ are anisotropic material parameters. The thermodynamic conjugate forces, namely the stress tensor and the energy release rate tensor, are obtained by taking the partial derivatives of the free energy potential with respect to the associated state variables: the strain tensor $\boldsymbol{\varepsilon}$ and the second order damage tensor $\boldsymbol{\omega}$.
Thus, the stress tensor is determined by
\begin{align}\label{eq:stress_tensor}
    \boldsymbol{\sigma} = \dfrac{\partial \phi(\boldsymbol{\varepsilon}, \boldsymbol{\omega})}{\partial \boldsymbol{\varepsilon}} =\lambda\;\mathrm{tr} (\boldsymbol{\varepsilon}) \boldsymbol{I}+ 2 \mu \boldsymbol{\varepsilon} + g \boldsymbol{\omega}
    + \alpha \left[\mathrm{tr}(\boldsymbol{\varepsilon} \cdot \boldsymbol{\omega} ) \boldsymbol{I} + \mathrm{tr}(\boldsymbol{\varepsilon})\boldsymbol{\omega} \right] + 2 \beta (\boldsymbol{\varepsilon} \cdot \boldsymbol{\omega} + \boldsymbol{\omega} \cdot \boldsymbol{\varepsilon}) \ .
\end{align}
The second-order energy release rate tensor is derived as
\begin{equation}\label{eq:energy_release_rate_tensor}
    \boldsymbol{Y} := -\dfrac{\partial \phi(\boldsymbol{\varepsilon}, \boldsymbol{\omega})}{\partial \boldsymbol{\omega}} = - g \boldsymbol{\varepsilon} - \alpha (\mathrm{tr}\boldsymbol{\varepsilon})\boldsymbol{\varepsilon} - 2 \beta (\boldsymbol{\varepsilon} \cdot \boldsymbol{\varepsilon}  )\ .
\end{equation}

To separate the compressive and tensile loading regimes, the strain tensor is decomposed into its positive component, $\boldsymbol{\varepsilon}^+$, and negative component, $\boldsymbol{\varepsilon}^-$. The damage evolution is specifically governed by the positive strain tensor $\boldsymbol{\varepsilon}^+$, capturing the tensile regime. The elastic domain is defined by a yield function $f$, which is expressed as
\begin{align}
\label{eq:yield_function}
f(\boldsymbol{\varepsilon},\boldsymbol{\omega} ) &= \sqrt{I^*_2} - k(\boldsymbol{\omega}) \leq 0 \ ,
\end{align}
where $I^*_2$ denotes the second invariant of the positive part of the energy release rate tensor $\boldsymbol{Y}^+$ and is defined as
\begin{equation}\label{eq:second_invariant}
I^*_2 = \dfrac{1}{2} \mathrm{tr}(\boldsymbol{Y}^+ \cdot \boldsymbol{Y}^+) = \dfrac{1}{2} \; g^2\; \mathrm{tr}(\boldsymbol{\varepsilon}^+ \cdot \boldsymbol{\varepsilon}^+) \ .
\end{equation}
The damage threshold function $k(\boldsymbol{\omega})$ is formulated as
\begin{equation}\label{eq:damage_threshold}
k(\boldsymbol{\omega}) = C_0 - C_1 \; tr(\boldsymbol{\omega}) \ ,
\end{equation}
$C_0$, and $C_1$ are material parameters. Therefore, the yield function can be expressed as
\begin{align}
\label{eq:yield_function_2}
f(\boldsymbol{\varepsilon},\boldsymbol{\omega} ) &= \dfrac{g}{\sqrt{2}} \tilde{\boldsymbol{\varepsilon}} - \left[ C_0 - C_1 \; tr(\boldsymbol{\omega}) \right] \ ,
\end{align}
where $\tilde{\varepsilon}$ is the equivalent strain proposed by Mazars \cite{mazars1989continuum}, defined as
\begin{equation}
\tilde{\varepsilon}= \sqrt{ (\langle \varepsilon_1 \rangle^+)^2 +(\langle \varepsilon_2 \rangle^+)^2 +(\langle \varepsilon_3 \rangle^+)^2} \ ,
\end{equation}
hereby, $\langle \varepsilon_1 \rangle^+$ represents the positive part of the eigenvalue of the strain tensor $\boldsymbol{\varepsilon}$. The extension of the model to account for fatigue loading follows the concept proposed by Marigo~\cite{marigo1985modelling}, which allows for the evolution of the damage parameter even within the damage yield surface. This approach has been widely adopted by various authors~\cite{Carrara_2020, KHALIL_2022, Baktheer_2024_PFM, Oneschkow_Loehnert_2022, Kristensen2023, KALINA_2023}. The yield limit is replaced by an irreversibility loading-unloading criterion given as
\begin{equation}
\begin{cases}
\text{if \quad $\dfrac{\partial f(\boldsymbol{\varepsilon},\boldsymbol{\omega} )}{\partial \boldsymbol{Y}^+}\cdot \boldsymbol{\dot{Y}}^+ > 0$:}  \Longrightarrow &\text{loading stage ($\boldsymbol{\dot{\omega}} > 0 $)}\\[5mm]
\text{if \quad $\dfrac{\partial f(\boldsymbol{\varepsilon},\boldsymbol{\omega} )}{\partial \boldsymbol{Y}^+}\cdot \boldsymbol{\dot{Y}}^+ \leq 0$:} \Longrightarrow  &\text{unloading stage ($\boldsymbol{\dot{\omega}} = 0 $)}.
\end{cases}
\end{equation}
Thus, damage progresses during loading and reloading without the need to reach a specified damage threshold. The damage evolution equation can be derived by differentiating Eq.~(\ref{eq:yield_function}) with respect to the energy release rate tensor as
\begin{equation}\label{eq:damage_1}
\boldsymbol{\dot{\omega}} = \dot{\lambda}_\mathrm{\omega} \; \dfrac{\partial f(\boldsymbol{\varepsilon},\boldsymbol{\omega} )}{\partial \boldsymbol{Y}^+} \ ,
\end{equation}
the damage multiplier $\dot{\lambda}_\mathrm{\omega}$ is determined by the consistency condition $\dot{f}(\boldsymbol{\varepsilon}, \boldsymbol{\omega}) = 0$, leading to the following equation
\begin{equation} \label{eq:multiplier}
\dot{\lambda}_\mathrm{\omega} = \frac{1}{h}\left(\frac{f}{K} \right)^n \left[ \frac{\partial f(\boldsymbol{\varepsilon},\boldsymbol{\omega} )}{\partial \boldsymbol{Y}^+}\cdot \boldsymbol{\dot{Y}}^+  \right] \ ,
\end{equation} 
where $h$, $n$, and $K$ are the material parameters. 
By combining the relations from equations (\ref{eq:damage_1}), (\ref{eq:multiplier}), and (\ref{eq:yield_function_2}), the resulting damage evolution equation is obtained as
\begin{equation}\label{eq:damage_2}
\boldsymbol{\dot{\omega}} = \left(\frac{f}{K} \right)^n \cdot \frac{\left[ \boldsymbol{\varepsilon}^+ : \boldsymbol{\dot{\varepsilon}}^+ \right]}{C_1 \;\mathrm{tr}(\boldsymbol{\varepsilon}^+) }  \cdot \frac{\boldsymbol{\varepsilon}^+}{\sqrt{2 \; \mathrm{tr}(\boldsymbol{\varepsilon}^+ \cdot \boldsymbol{\varepsilon}^+)}}
\end{equation}

\subsection{Dimensional reduction to uniaxial stress state in fatigue modeling}

To model fatigue behavior under compressive loading, the problem is simplified by assuming a uniform stress and strain distribution within the specimen. This leads to the assumption that $\varepsilon_{1} < 0$ (indicating compressive strain along the loading axis) and $\varepsilon_{2} = \varepsilon_{3} > 0$ (indicating the strains in the transverse directions are equal and positive). Since damage is governed solely by the positive part of the strain tensor, the stress, strain, and damage tensors can then be expressed as follows
\begin{equation} \label{eq:stress_strain_tensor}
\boldsymbol{\sigma} = \begin{bmatrix}
\sigma_{1} & 0 & 0 \\
0 & 0 & 0 \\
0 & 0 & 0
\end{bmatrix},\quad \quad
\; \boldsymbol{\varepsilon} = \begin{bmatrix}
\varepsilon_{1} & 0 & 0 \\
0 & \varepsilon_{2} & 0 \\
0 & 0 & \varepsilon_{3}
\end{bmatrix},\quad \quad
\; \boldsymbol{\omega} = \begin{bmatrix}
0 & 0 & 0 \\
0 & \omega_{2} & 0 \\
0 & 0 & \omega_{2}
\end{bmatrix}.
\end{equation} 
Based on the constitutive equation Eq.~(\ref{eq:stress_tensor}), the stress components in Eq.~(\ref{eq:stress_strain_tensor}) can be rewritten as
\begin{align} 
&\sigma_1 = (\lambda +2 \mu) \varepsilon_1 + 2(\lambda+ \alpha \; \omega_2) \varepsilon_2 \label{eq:sigma_1}\\
&\sigma_2 = 0 = \left[2(\lambda + \mu)  + 4(\alpha +\beta)\omega_2 \right]\varepsilon_1 + (\lambda+ \alpha \; \omega_2)\varepsilon_1 + g \;\omega_2  \label{eq:sigma_2} 
\end{align} 
From Eq.~(\ref{eq:sigma_1}) the strain component $\varepsilon_2$ can be resolved as follows
\begin{equation} \label{eq:epsilon_1}
\varepsilon_2 = \frac{(\lambda+ \alpha \; \omega_2)\sigma_1 + g \;\omega_2 (\lambda +2 \mu)}{(\lambda +2 \mu) \left[2 (\lambda + \mu) + 4 (\alpha+\beta)\omega_2 \right] - 2(\lambda+\alpha \omega_2)^2} 
\end{equation}
The yield function can then be reduced to
\begin{equation}
f(\boldsymbol{\varepsilon},\boldsymbol{\omega} ) = |g|\;\varepsilon_2 - (C_0 + 2 C_1 \; \omega_2) \ .
\end{equation}
Thus, the damage evolution law can be expressed as
\begin{equation}
\dot{\omega}_2 = \frac{|g|}{2 C_1} \left(\frac{f}{K}\right)^n \left\langle - \frac{\lambda + \alpha \;\omega_2}{\kappa} \dot{\sigma_1} \right\rangle^+ ,
\end{equation}
where $\kappa$ is defined as
\begin{align}
\kappa = (\lambda + 2 \mu)  2(\lambda +  \mu)  + 4(\alpha +\beta)\omega_2 - \alpha \frac{g}{2C_1}(2\varepsilon_2 + \varepsilon_1) 
 - \frac{g^2}{2C_1}  -2(\lambda + \alpha \; \omega_2)^2
\end{align}

\subsection{Numerical Implementation}
The model is prepared for integration into a standard time-stepping algorithm that captures the state of a compressive specimen using two state variables and a single degree of freedom. This enables efficient, cycle-by-cycle simulations of high-cycle fatigue loading, which can be executed in a matter of minutes on standard computational platforms. Due to its incremental formulation, the model is versatile and can accommodate various loading scenarios.
The dimensional reduction to a uniaxial stress state simplifies the complexity of multi-dimensional stress conditions, while still preserving the material's key behavior under uniaxial loading. This approach not only enhances the efficiency of simulations but also facilitates the generation of data suitable for training ML models.

\begin{figure*}[!b]
\centerline{
\includegraphics[width=1.0\textwidth]{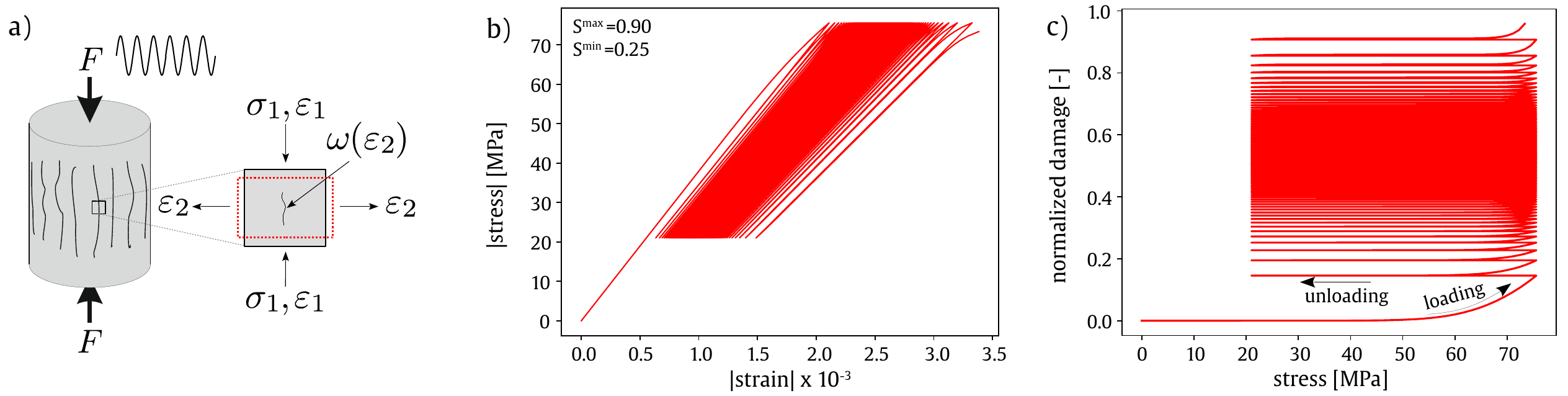}}
\caption{Elementary example of a fatigue simulation using the applied anisotropic damage fatigue model: a) uniaxial compression test setup and the used anisotropic modeling approach;\; b) stress-strain curve;\; c) corresponding damage evolution during the loading history}
\label{f:Alliche_model_example}
\end{figure*}

To investigate the fundamental behavior obtained by the described fatigue model, an example of uniaxial fatigue loading on a concrete cylinder with a compressive strength of $f_{\mathrm{c}}=84$~MPa~\cite{KIM19961513} is presented in Fig.~\ref{f:Alliche_model_example}. The simulation considers a fatigue loading scenario with a maximum stress range 
$S^\mathrm{max}=0.9$ and minimum stress range $S^\mathrm{min}=0.25$, reflecting typical loading conditions. The corresponding stress-strain response over the complete loading history, along with the evolution of damage, are shown in Figs.~\ref{f:Alliche_model_example}b and \ref{f:Alliche_model_example}c, respectively. These figures illustrate how the material behaves under cyclic loading, showing key characteristics such as elastic deformation, damage accumulation, and the transition to failure under sustained cyclic loading. The damage evolution reflects the underlying mechanics of fatigue deterioration in concrete, providing insight into how the material degrades over multiple cycles.

\begin{figure*}[!b]
\centerline{
\includegraphics[width=1.0\textwidth]{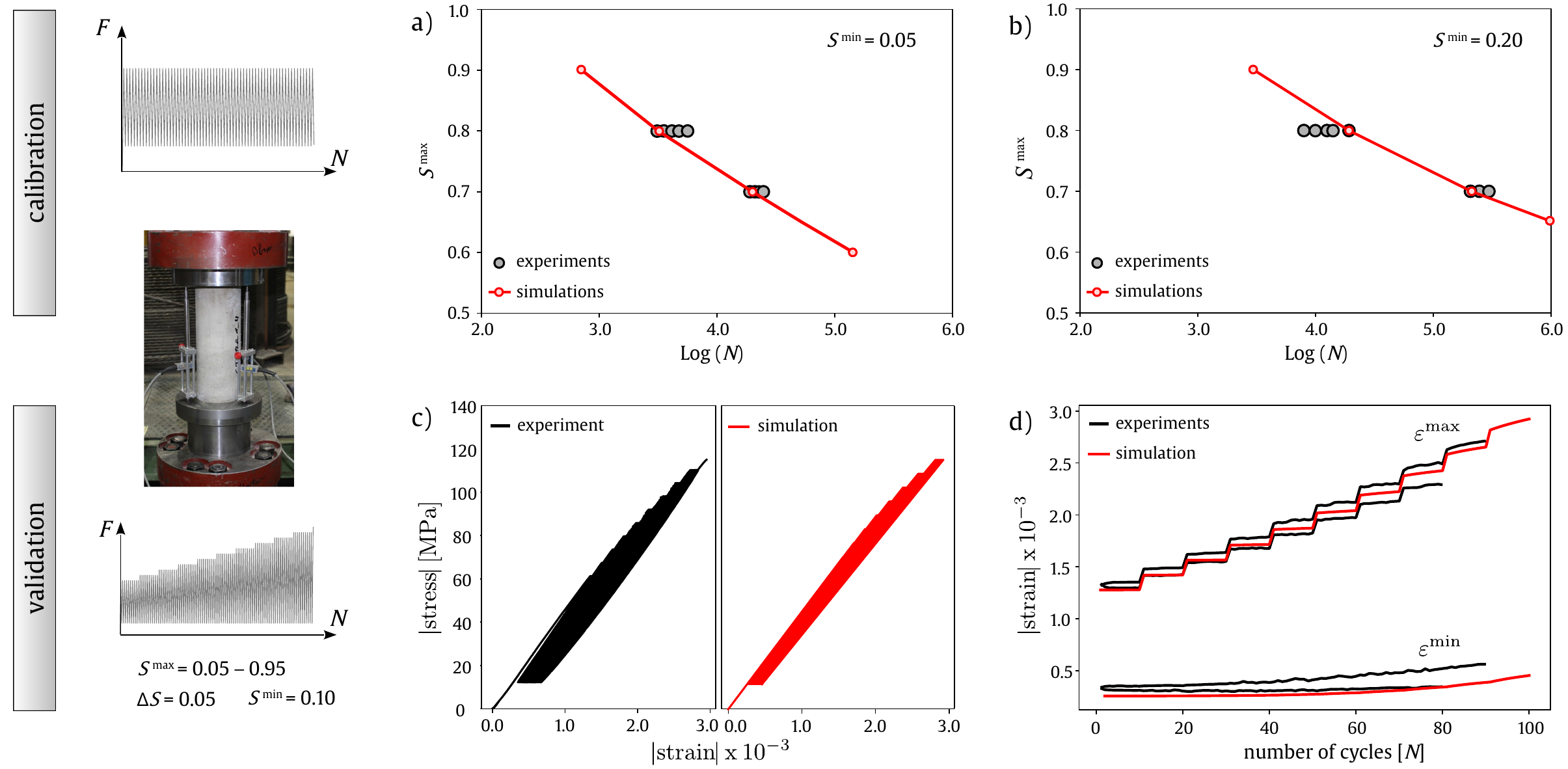}}
\caption{Calibration and validation procedure of the applied anisotropic damage fatigue model: a, b) calibrated S-N curves representing the fatigue behavior under constant amplitude loading scenarios;\; c) stress-strain curves: comparison between experimental data and model predictions; and d) corresponding fatigue creep curves for upper and lower load levels
}
\label{f:Calibration_validation}
\end{figure*}

\subsection{Model calibration and validation}

The model parameters are calibrated based on the experimental program reported in~\cite{schneider2018untersuchungen}. This study investigated the compressive fatigue behavior of concrete class C100/115 under constant amplitudes. The tests were conducted on cylindrical specimens with dimensions of $100 \times 300$~mm at a loading frequency of 5~Hz.
The Lam$\Acute{\mathrm{e}}$ constants, $\lambda$ and $\mu$, are derived from the elastic properties: Young's modulus $E = 49000$~MPa and Poisson’s ratio $\nu = 0.2$. The remaining six parameters are identified through curve fitting of the fatigue test results, as illustrated in Figs.~\ref{f:Calibration_validation}a and \ref{f:Calibration_validation}b. 
The model is compared to experimental S-N curves at lower loading levels of 0.05 and 0.2, demonstrating its predictive accuracy. The calibrated material parameters are summarized in Table~\ref{t:alliche_model_parameters}.

The model is validated using experimental data from a test program with stepwise increasing maximum load levels, as described in~\cite{Baktheer_2022_tip_bearing} and illustrated in Fig.~\ref{f:Calibration_validation}. The predicted response is presented as a red stress-strain curve in Fig.~\ref{f:Calibration_validation}c, alongside experimental results plotted as black curves. The corresponding fatigue creep curves for the upper and lower load levels are shown in Fig.~\ref{f:Calibration_validation}d. The comparison indicates a reasonable agreement between the numerical predictions of the used anisotropic fatigue model and the experimental results.

\begin{table*}[!h]
\renewcommand{\arraystretch}{1.3} 
 \begin{center}
 \caption{Calibrated material parameters of the used anisotropic fatigue model for concrete class C120}\label{t:alliche_model_parameters}
 \begin{tabular}{ m{2.4cm} m{1.3cm} m{1.3cm} m{1.3cm} m{1.5cm} m{1.0cm} m{1.3cm} m{1.3cm} m{1.3cm} m{0.6cm} }
    \toprule
    \textbf{Parameter} & $\lambda$&$\mu$ & $g$ & $K$ &$C_0$& $C_1$& $\alpha$& $\beta$ &$n$\\ \hline
   \textbf{Value} &  12500 & 18750 & -10.0 & 0.00485 & 0.0 & 0.0019 & 2237.5 & -2116.5 & 10\\
  \bottomrule
 \end{tabular}
 \end{center}
\end{table*}

\subsection{Effect of fatigue loading sequence }
\begin{figure*}[!t]
\centerline{
\includegraphics[width=1.0\textwidth]{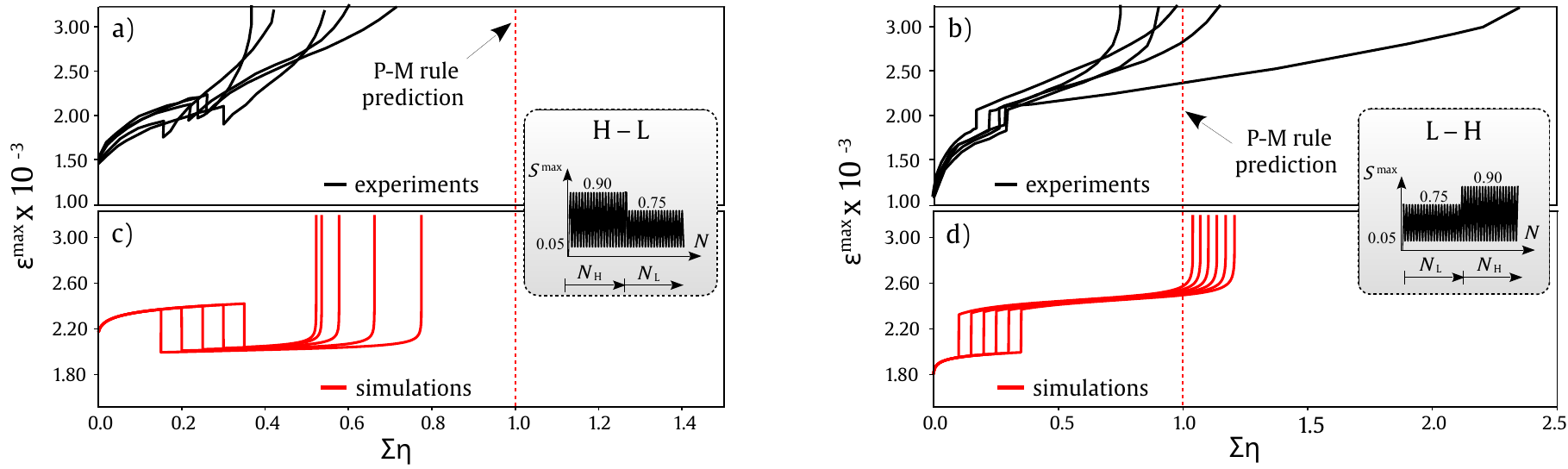}}
\caption{Qualitative evaluation of the loading sequence effects on fatigue behavior, comparing numerical simulations and experimental data from~\cite{holmen1982fatigue}: a) experimental results for H-L loading sequence;\; b) experimental results for L-H loading sequence;\; c) numerical results for H-L loading sequence; and d) numerical results for L-H loading sequence}
\label{f:Qualitative_comparison}
\end{figure*}

The experimental findings presented in~\cite{holmen1982fatigue}, investigating the effect of loading sequence on the fatigue behavior of concrete class C40, are presented in Figs.~\ref{f:Qualitative_comparison}a and \ref{f:Qualitative_comparison}b. 
These results are depicted as fatigue creep curves for the used  (H-L) and (L-H) loading scenarios. The consumed fatigue lifetime in the first loading range was set to specific values (15\%, 20\%, 25\%, 27\%, and 30\%), and the lifetime on the horizontal axis is normalized with respect to predictions made using the Palmgren-Miner (P-M) rule.
According to this rule, fatigue failure occurs when the fatigue loading history, comprising $n$ distinct loading ranges, satisfies the following condition:
\begin{equation} \label{eq:P-M_rule}
   \sum \eta =  \sum_{\mathrm{i}=1}^n \dfrac{N_{\mathrm{i}}}{N^{\mathrm{f}}_{\mathrm{i}}} = 1,
\end{equation}
where $N_\mathrm{i}$ represents the number of applied cycles for the $\mathrm{i}$-th loading range, and $N^\mathrm{f}_\mathrm{i}$ denotes the number of cycles to failure under constant amplitude loading for that range. Thus, the ratio $N_{\mathrm{i}} /N^{\mathrm{f}}_{\mathrm{i}}$ indicates the consumed fatigue lifetime for the corresponding loading range.

The (H-L) loading scenario (Fig.~\ref{f:Qualitative_comparison}a) resulted in significantly shorter fatigue lifetime compared to the P-M rule predictions. Conversely, the (L-H) scenario yielded fatigue lives that, in some cases, exceeded the P-M rule estimates (Fig.~\ref{f:Qualitative_comparison}b).
The numerical model predictions for concrete class C110/115, using the parameters provided in Table~\ref{t:alliche_model_parameters}, under similar loading scenarios are shown in Figs.~\ref{f:Qualitative_comparison}c and \ref{f:Qualitative_comparison}d. Despite the difference in concrete classes between the experiments (C40) and the simulations (C110/115), the model qualitatively reproduce the experimentally observed trends: a relative reduction in fatigue lifetime for (H-L) loading and an extension for (L-H) loading. These trends were experimentally confirmed in recent test programs presented in~\cite{Baktheer_2021_3, BECKS_2024}.

While the validity of the numerical model cannot be conclusively established based on a single experimental dataset with differing concrete classes,  it nonetheless demonstrates a reasonable ability to capture the observed qualitative trends. In light of the scarcity of experimental data, the model offers a valuable tool for systematically studying the loading sequence effect across a wide range of configurations. This can serve as a basis for the training of a neural network that enables the prediction of remaining fatigue lifetime under general loading scenarios, serving as fatigue damage accumulation assessment rule for concrete.

\section{Data driven and physics-based neural networks}
\label{sec:neural_networks}

%at a single material point

This section outlines the development of machine learning models for predicting the remaining fatigue lifetime of structures, accounting for the nonlinear effects of the loading sequence.
As highlighted in the introduction, the limited availability of experimental data on the effects of loading sequences, combined with the broad parameter space of complex loading scenarios, poses significant challenges for deriving reliable damage accumulation rules for concrete fatigue. Existing approaches in the literature, such as those in~\cite{grzybowski1993damage, shah1984predictions}, are often oversimplified and insufficient for capturing the intricacies of concrete fatigue behavior~\cite{BAKTHEER2019}.
Machine learning models, particularly physics-based machine learning approaches, offer a flexible and powerful framework for developing predictive models capable of addressing these complexities.
To address the issue of data scarcity from experimental studies, the calibrated and validated anisotropic fatigue model described in Sec.~\ref{sec:fatigue_model} is employed to generate a comprehensive dataset for training the machine learning models.

\begin{figure*}[!b]
\centerline{
\includegraphics[width=1.0\textwidth]{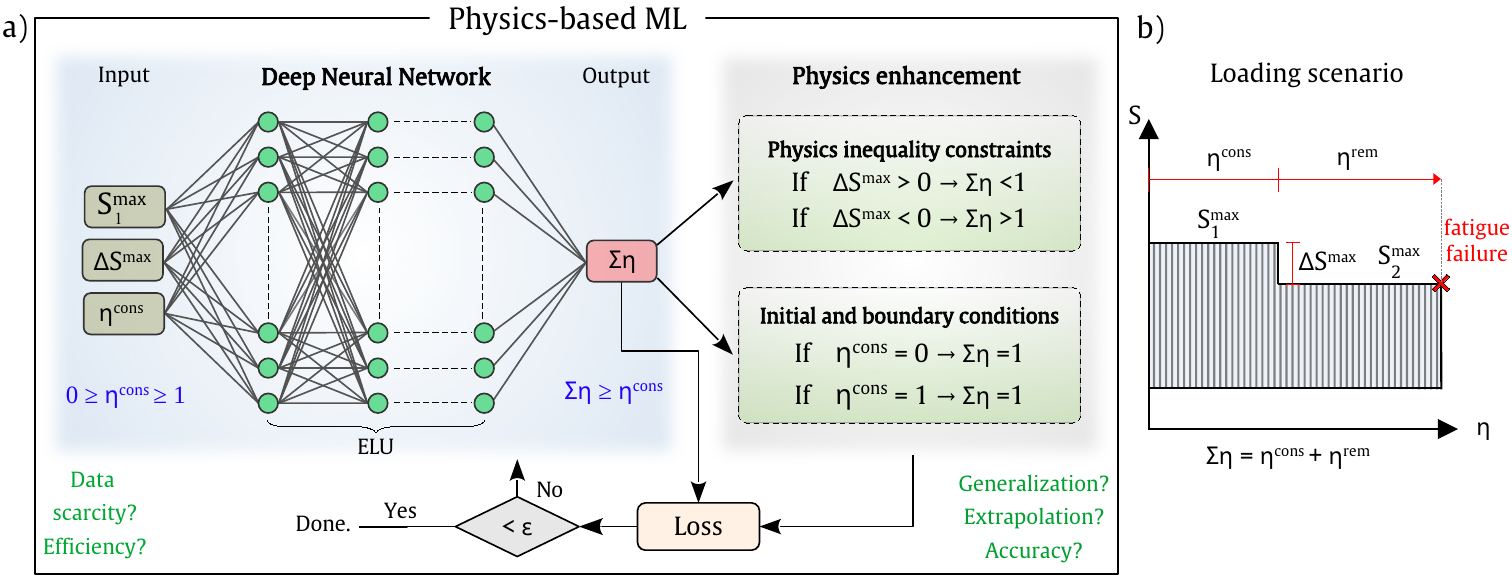}}
\caption{Schematic representation of the developed physics-based machine learning model for predicting fatigue lifetime: a) model architecture for $\phi$ML approach, incorporating a customized loss function with imposed physical constraints;\;
b) characteristics of a two-level loading scenario with a single load jump
}
\label{f:phyNN_architecture}
\end{figure*}

% and forward propagation

\subsection{Neural network architecture}
\label{sec:NN_architecture}

To achieve accurate predictions of the remaining fatigue lifetime under non-uniform loading scenarios, we employ a physics-based machine learning ($\phi$ML) framework.
This approach integrates data-driven deep learning with physics-based constraints, ensuring that the predictions align with established principles and empirical knowledge of the fatigue behavior of structures, particularly under complex non-uniform loading conditions. By embedding physics-based constraints into the learning process, the $\phi$ML can reduce the risk of overfitting while maintaining consistency with underlying mechanical phenomena.
The model is designed as a multilayer perceptron (MLP) to predict the remaining fatigue lifetime \( \eta^{\text{rem}} \) for a given load scenario with two different loading ranges (see Fig.~\ref{f:phyNN_architecture}b), with a focus on accurately capturing fatigue behavior while considering the effects of the load sequence.

\bigskip
The input features represent key parameters that define the loading scenario (Fig.~\ref{f:phyNN_architecture}b) and include the following:
\begin{itemize}
    \item The maximum stress ratio of the first loading range, \( S^{\text{max}}_1 \).
    \item The magnitude of the load jump between the two loading ranges, \( \Delta S^{\max} = S^{\text{max}}_1 - S^{\text{max}}_2 \).
    \item The consumed fatigue lifetime prior to the load jump, \( \eta^\text{cons} \) (see Fig.~\ref{f:phyNN_architecture}).
\end{itemize}
It is important to note that the presented example of the $\phi$ML model focuses on loading scenarios where the lower load level is kept constant, i.e., $\Delta S^{\min} = 0$, in order to simplify the conceptual framework. However, the proposed approach can be readily extended to accommodate loading scenarios with $\Delta S^{\min} \neq 0$, demonstrating its flexibility and adaptability to more complex conditions.

The framework is designed to accurately represent fatigue behavior, particularly by incorporating the effects of loading sequences. This approach aims to enable the $\phi$ML to account for non-linear damage accumulation effects, which are typically difficult to model using conventional empirical methods. By integrating physics-guided constraints, as will be discussed later, the framework is intended to address scenarios with limited experimental data, ultimately enhancing the robustness and generalizability of its predictions.
%

%#### Neural Network Model Definition

The feed-forward neural network (FFNN) model for both the purely data-driven ML and physics-based $\phi$ML is defined mathematically as follows:
\begin{equation}
    FFNN: \quad (S^{\text{max}}_1, \Delta S^{\max}, \eta^\text{cons})     \quad \longmapsto \quad  
 \sum \eta =  \eta^\text{cons} + \eta^{\text{rem}}
\end{equation}
the neural network takes a three-dimensional input vector ($S^{\text{max}}_1, \Delta S^{\max}, \eta^\text{cons}$) and outputs a scalar value. The model's output,  $y_{\text{pred}} = \sum \eta$, is the predicted cumulative fatigue life, computed as the sum of the consumed fatigue lifetime $\eta^\text{cons}$ and the remaining fatigue lifetime $\eta^\text{rem}$.
The FFNN employs a fully connected deep neural network with 10 hidden layers, each consisting of 16 neurons, and utilizing the Exponential Linear Unit (ELU) activation function.
The ELU activation function is mathematically expressed as:
\begin{equation}
  \text{ELU}(x) = 
  \begin{cases} 
      x & \text{if } x > 0 \\
      \zeta (e^x - 1) & \text{if } x \leq 0, 
  \end{cases}
\end{equation}
where the constant $\zeta$ is set to 1.0. The ELU activation function is particularly advantageous for capturing complex patterns, as it mitigates issues like vanishing gradients and enhances the network’s ability to model highly nonlinear relationships. It should be noted alternative activation functions, such as the Rectified Linear Unit (ReLU), can be equally suitable for the task addressed in this study.

During forward propagation, the input vector \(\mathbf{x} = [S^{\text{max}}_1, \Delta S^{\max}, \eta^\text{cons}]\) sequentially passes through each layer of the network. For a given hidden layer \( l \), the computation at each step is:  
\begin{equation}
\mathbf{z}^{(l)} = \text{ELU}(\mathbf{W}^{(l)} \mathbf{a}^{(l-1)} + \mathbf{b}^{(l)}), 
\end{equation}  
where:  
\begin{itemize}
    \item \( \mathbf{W}^{(l)} \) is the weight matrix for layer \( l \),  
    \item \( \mathbf{b}^{(l)} \) is the bias vector for layer \( l \),  
    \item \( \mathbf{a}^{(l-1)} \) is the activation output of the previous layer (or the input vector \(\mathbf{x}\) for the first layer),  
    \item \( \mathbf{z}^{(l)} \) is the output of layer \( l \) after applying the ELU activation.  
\end{itemize}
This computation propagates through all hidden layers, ultimately reaching the output layer, where a linear transformation is applied to produce the final scalar output \( y_{\text{pred}} \).  The architecture is designed to capture highly complex, non-linear relationships between input parameters and accumulated fatigue damage. The FFNN structure, combined with ELU activation function, equips the network with the flexibility and capacity needed to handle the complex effects of loading sequence and damage accumulation on concrete fatigue behavior. Additionally, the physics-based framework of the $\phi$ML integrates domain knowledge into the learning process, as will be elaborated in the next subsection.

\begin{figure*}[!b]
\centerline{
\includegraphics[width=1.0\textwidth]{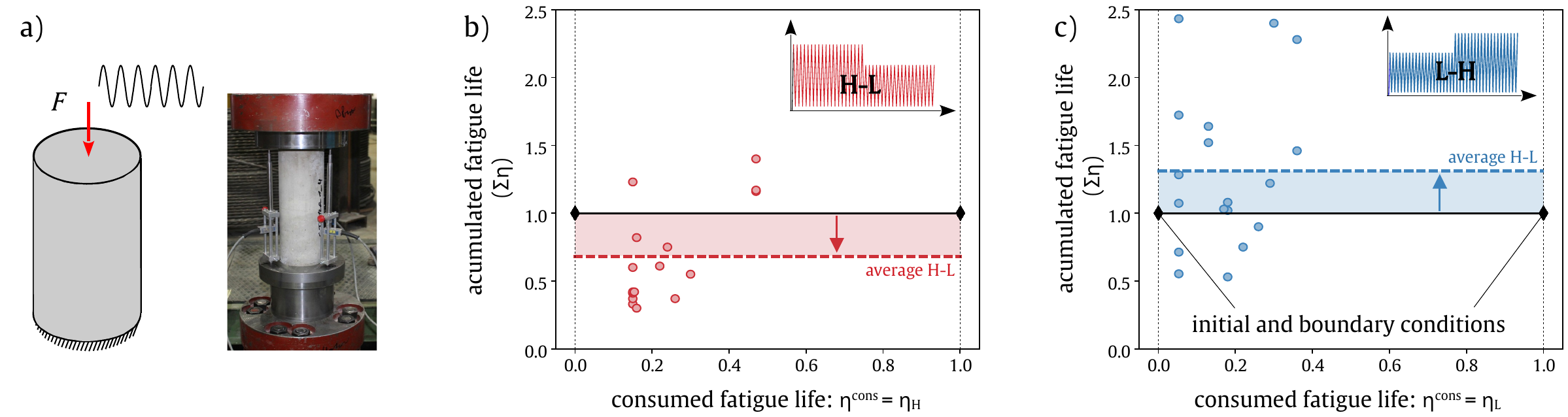}}
\caption{
Experimentally observed trends for fatigue lifetime~\cite{Baktheer_2021_3} imposed as physical constraints within the developed machine learning model: a) test setup for concrete fatigue under compression;\; b) H-L scenario;\; c) L-H scenario}
\label{f:pyhsical_constraints}
\end{figure*}

\subsection{Physics-based constraints}
To enforce physical knowledge on loading sequence effects in fatigue lifetime, physics-based constraints are embedded throughout the model architecture, including a custom loss function. Experimental evidence~\cite{Baktheer_2021_3}, published recently by the authors, has demonstrated distinct trends in fatigue behavior for different loading sequences, such as high-low (H-L) and low-high (L-H) loading scenarios. These trends are illustrated in Fig.~\ref{f:pyhsical_constraints}, where H-L sequences result in a shortened fatigue lifetime ($\sum \eta < 1$) and L-H sequences exhibit an extended lifetime ($\sum \eta > 1$), compared to the linear damage accumulation predicted by the Palmgren-Miner (P-M) rule ($\sum \eta = 1$), as depicted in Figs.~\ref{f:pyhsical_constraints}b and~\ref{f:pyhsical_constraints}c, respectively. Similar experimental observations have been reported in~\cite{holmen1982fatigue, BECKS_2024, jinawath1974cumulative, hoff1984testing, petkovic1990fatigue}.  
This experimentally observed behavior is implemented as physical constraints within the machine learning (ML) model architecture, ensuring that the training process is guided by known fatigue characteristics. Additionally, the boundary conditions for consumed fatigue life, $\eta^\mathrm{cons} = 0$ (no fatigue lifetime consumed in the first range) and $\eta^\mathrm{cons} = 1$ (complete fatigue lifetime consumed in the second loading range) highlighted in Fig.~\ref{f:phyNN_architecture}, are incorporated into the developed $\phi$ML model architecture.  
To ensure consistency of the $\phi$ML model, physical constraint on the output value is considered to fulfill the condition $\sum \eta \geq \eta^{\mathrm{cons}}$ as highlighted in Fig.~\ref{f:phyNN_architecture}a.
Building on this physics-based model architecture, the training process is further guided by a custom loss function, $\mathcal{L}_{\text{cust}}$, which combines a conventional mean squared error (MSE) term with physics-based constraints and boundary conditions, as outlined in Fig.~\ref{f:phyNN_architecture}a. 
This design penalizes deviations from physically plausible fatigue behavior and ensures the model remains consistent with experimental trends. The loss function incorporates the following key components:

\paragraph{\textbf{Data loss (Mean Squared Error)}} 
The primary component of the loss function is the data loss, calculated as the mean squared error (MSE) between the predicted fatigue lifetime values ($y_{\text{pred}}$) and the ground truth values ($y_{\text{target}}$). This term ensures that the model accurately fits the training dataset, which contains $N$ data points:
\begin{equation}
   \mathcal{L}_{\text{data}} = \frac{1}{N} \sum_{i=1}^{N} (y_{\text{pred}, i} - y_{\text{target}, i})^2.
\end{equation}

\paragraph{\textbf{Physics-based constraints loss}}
The experimental trends illustrated in Fig.~\ref{f:pyhsical_constraints} are incorporated into the loss function as inequality constraints. These constraints enforce the expected behavior of cumulative fatigue lifetime under varying loading scenarios, expressed as:

\begin{equation}
    \begin{cases} 
      \sum \eta < 1, & \text{if } \;\; \Delta S^\mathrm{max} > 0, \\[2mm]
      \sum \eta > 1, & \text{if } \;\; \Delta S^\mathrm{max} < 0. 
    \end{cases}
\end{equation}

The physics-based constraints guide the model to capture realistic fatigue accumulation behavior. Specifically:
\begin{itemize}
    \item Positive constraint: For decreasing \( S^{\text{max}} \) (i.e. \( \Delta S^{\text{max}} > 0\)),  cumulative fatigue lifetime predictions should approach values below 1, reflecting accelerated fatigue accumulation.
    \item Negative constraint: For increasing \( S^{\text{max}} \) (i.e. \( \Delta S^{\text{max}} < 0\)), cumulative fatigue lifetime predictions should exceed 1, reflecting decelerated fatigue accumulation.
\end{itemize}

The physics-based constraint term, \( \mathcal{L}_{\text{const}} \), is defined as:

\begin{equation}
   \mathcal{L}_{\text{const}} = w_{\text{const}} \left( \sum_{\text{pos}} (1 - y_{\text{pred}})^2 + \sum_{\text{neg}} (y_{\text{pred}} - 1)^2 \right),
\end{equation}
where \( w_{\text{const}} \) is a weight parameter controlling the influence of the constraints. The first term penalizes deviations from the expected cumulative fatigue lifetime for decreasing \( S^{\text{max}} \), while the second term ensures accurate predictions for increasing \( S^{\text{max}} \).

\paragraph{\textbf{Initial and boundary conditions based loss}}
To ensure that the model’s predictions align with physical boundary behaviors, the following boundary conditions are considered:
\begin{equation}\label{boundary-conditions}
    \begin{cases} 
    \sum \eta = 1, & \text{if }\;\; \eta^\mathrm{cons} = 0,\\[2mm]
      \sum \eta = 1, & \text{if }\;\; \eta^\mathrm{cons} = 1. 
    \end{cases}
\end{equation}
The initial condition for $\eta^\mathrm{cons} = 0$ represents that no fatigue lifetime has been consumed in the first loading range, then the cumulative lifetime should be close to 1, reflecting the initial full fatigue lifetime (i.e. $\sum \eta = \eta^\mathrm{rem} = 1$).
The boundary condition $\eta^\mathrm{cons} = 1$ represents that all fatigue lifetime is consumed at the first load range, then the model should predict a cumulative lifetime of approximately 1, implying negligible remaining fatigue lifetime ($\eta^\mathrm{rem} = 0$). The boundary loss term can be then written as:
\begin{equation}
   \mathcal{L}_{\text{bound}} = w_{\text{bound}} \left( (1 - y_{\text{pred}} |_{\eta^\mathrm{cons} = 0})^2 + (y_{\text{pred}} |_{\eta^\mathrm{cons} = 1} - 1)^2 \right)
\end{equation}

To further enhance the predictive performance of the $\phi$ML model, an algorithmic enhancement termed Sparse-Domain Regularization is incorporated and applied to regions with sparse data coverage.  
This regularization complements the physical boundary constraints imposed at {\it feature 3} ($\eta^\mathrm{cons}$) by guiding the network toward reasonable predictions in sparse regions, reducing the likelihood of unrealistic jumps or irregular patterns in low-density areas. Specifically, it targets the domain where $\eta^\mathrm{cons} < 0.35$ (see Fig.~\ref{f:data_generation}c).  
The combination of strict boundary constraints $\eta^\mathrm{cons} = 0$ in (\ref{boundary-conditions}) and regularization within the low-density region ensures smooth interpolation between these regions.  

The loss term for sparse-domain regularization is defined as:  
\begin{equation}
   \mathcal{L}_{\text{spars}} = w_{\text{spars}} \sum_{\text{spars}} (1 - y_{\text{pred}})^2
\end{equation}  

The overall custom loss function combines all the relevant components as follows:  
\begin{equation}
\mathcal{L}_{\text{total}} = \mathcal{L}_{\text{data}} + \mathcal{L}_{\text{const}} + \mathcal{L}_{\text{bound}} + \mathcal{L}_{\text{spars}}
\end{equation}  
where $w_{\text{const}}$, $w_{\text{bound}}$, and $w_{\text{spars}}$ are hyperparameters adjusted to balance the influence of each term.

\subsection{Data generation for training the machine learning models}
\label{sec:data_generation}

The applied fatigue model, calibrated and validated as shown in Fig.~\ref{f:Calibration_validation}, is used to generate a dataset for training the machine learning models. To capture the effects of loading sequence on fatigue lifetime, loading scenarios involving two load levels (i.e., H-L or L-H) are constructed, ensuring a broad range of possible variations in the input features: $S^{\text{max}}_1$, $\Delta S^{\text{max}}$, and $\eta^{\text{cons}}$. The upper load limit $S^{\text{max}}_1$ is discretized into the following six levels:
\[
S^{\text{max}}_1 = \{0.65,\;  0.70,\;   0.75,\; 0.80,\; 0.85,\; 0.90\}
\]
These values represent subcritical loading levels, where fatigue damage occurs over a spectrum ranging from low-cycle to high-cycle fatigue as reported experimentally by many authors e.g.,~\cite{schneider2018untersuchungen, Baktheer_2021_3}.
The upper limit of the second load level $S^{\text{max}}_2$ (Fig.~\ref{f:phyNN_architecture}b) also varies between the same values as $S^{\text{max}}_1$.

For all generated scenarios, the lower load limit was fixed at $S^{\text{min}} = 0.2$. The fatigue lifetime $N^{\text{f}}$ for each loading level was determined based on the S-N curves presented in Fig.~\ref{f:Calibration_validation}b. To define the loading jump $\Delta S^{\text{max}}$ between two load levels, all possible combinations of $S^{\text{max}}_1$ values were systematically considered. The discrete values of $\Delta S^{\text{max}}$ for each $S^{\text{max}}_1$ are summarized in the following matrix:

\[
\begin{array}{c|ccccc}
S^{\text{max}}_1 & \multicolumn{5}{c}{\Delta S^{\text{max}}} \\[1mm] \hline
0.65 & \cellcolor{gray!30} -0.05 & \cellcolor{gray!30} -0.10 & \cellcolor{gray!30} -0.15 & \cellcolor{gray!30} -0.20 & \cellcolor{gray!30} -0.25   \\ 
0.70 & \cellcolor{white} 0.05  & \cellcolor{gray!30} -0.05  & \cellcolor{gray!30} -0.10  & \cellcolor{gray!30} -0.15  & \cellcolor{gray!30} -0.20    \\ 
0.75 & \cellcolor{white} 0.10  & \cellcolor{white} 0.05  & \cellcolor{gray!30} -0.05  & \cellcolor{gray!30} -0.10  & \cellcolor{gray!30} -0.15    \\ 
0.80 & \cellcolor{white} 0.15  & \cellcolor{white} 0.10  & \cellcolor{white} 0.05  & \cellcolor{gray!30} -0.05  & \cellcolor{gray!30} -0.10    \\ 
0.85 & \cellcolor{white} 0.20  & \cellcolor{white} 0.15  & \cellcolor{white} 0.10  & \cellcolor{white} 0.05  & \cellcolor{gray!30} -0.05   \\ 
0.90 & \cellcolor{white} 0.25  & \cellcolor{white} 0.20  & \cellcolor{white} 0.15  & \cellcolor{white} 0.10  & \cellcolor{white} 0.05    \\ 
\end{array}
\]
whereby positive values of $\Delta S^{\text{max}}$ correspond to H-L scenarios and negative values (highlighted in gray) to L-H scenarios. To further account for sequence effects, the number of loading cycles applied at the first loading level $S^{\text{max}}_1$ is varied to represent different fractions of consumed fatigue life. The consumed fatigue lifetime ratio $\eta^{\text{cons}}$ is defined as:
\[
\eta^{\text{cons}} = \frac{N_1}{N_1^\mathrm{f}}, \quad \eta^{\text{cons}} = \{ 0.0,\; 0.05,\; 0.10,\; \ldots,\; 1.0 \}, \quad \Delta \eta^{\text{cons}} = 0.05
\]
where the special cases of $\eta^{\text{cons}} = 0$ and $\eta^{\text{cons}} = 1$ represent scenarios in which the entire fatigue lifetime is consumed in a single load range.

Considering all possible combinations of the three features, a total of \textbf{630 loading scenarios} are generated. Each scenario underwent fatigue simulations using a cycle-by-cycle approach based on the fatigue model, allowing the determination of the remaining fatigue lifetime at the second loading range. The final output for each scenario is computed as the sum of the consumed fatigue lifetime in the first range and the remaining fatigue lifetime in the second range (see Fig.~\ref{f:phyNN_architecture}b). 

\begin{figure*}[!t]
\centerline{
\includegraphics[width=1.0\textwidth]{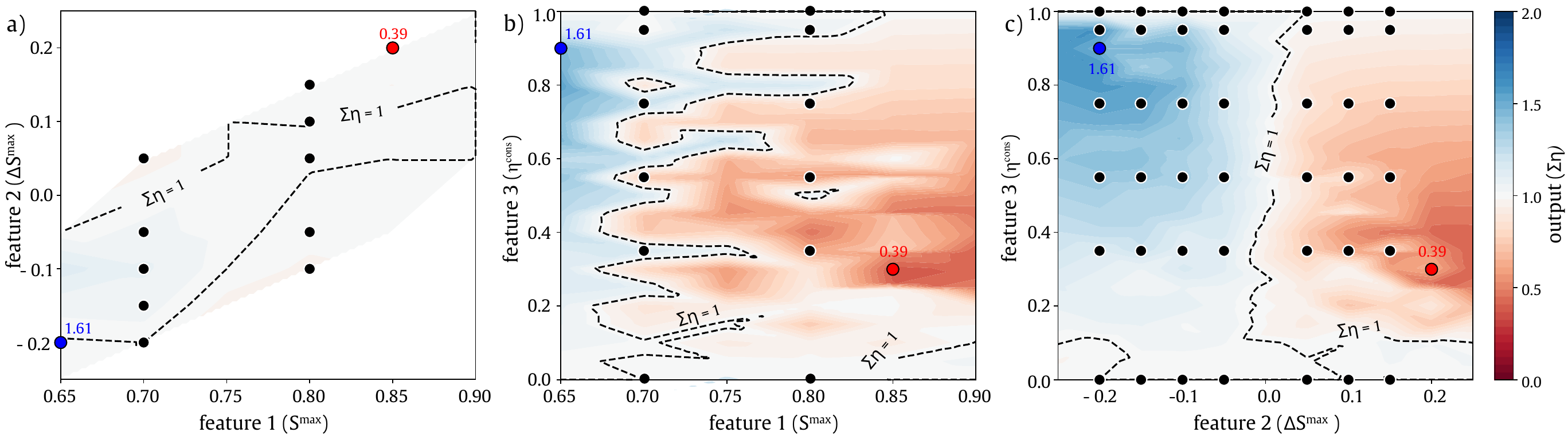}}
\caption{
Contour plots illustrating the relationship between the three input features and the output of the machine learning model, including interpolations between feature data points. These visualizations are based on data generated by the anisotropic fatigue model used to train the machine learning model: a) interaction between feature~1 and feature~2;\; b) interaction between feature~1 and feature~3;\; c) interaction between feature~2 and feature~3}
\label{f:data_generation}
\end{figure*}

To visualize the relationships between the three input features and the output values, contour plots of the generated data are presented in Fig.~\ref{f:data_generation}. The color scale in these plots indicates the deviation from Palmgren-Miner (P-M) rule predictions at $\sum \eta = 1$, with maximum and minimum values of 1.61 and 0.39, respectively, as shown by the colored circles in Fig.~\ref{f:data_generation}. These results highlight the highly nonlinear and complex influence of loading sequence effects on fatigue behavior, an aspect that remains insufficiently characterized in current concrete fatigue research.

Investigating the potential of physics-based machine learning models for fatigue lifetime prediction, particularly in cases with limited experimental data, involved preparing two datasets for model training:
\begin{itemize}
    \item \textbf{Large dataset:} This dataset includes the full \textbf{630 samples}, divided into training (70\%), validation (15\%), and testing (15\%). It provides a comprehensive basis for developing and evaluating data-driven models.
    \item \textbf{Small dataset:} A subset of \textbf{60 samples} is selected from the total dataset to examine the feasibility of training models with minimal data, mimicking scenarios where fatigue data is costly and scarce. These samples were chosen based on discrete feature values: specifically, $S^{\text{max}} = \{0.70, 0.80\}$, the corresponding values of $\Delta S^{\text{max}}_1$ and $\eta^{\text{cons}} = \{0.0, 0.35, 0.55, 0.75, 0.95, 1.0\}$, as highlighted by black circles in Fig.~\ref{f:data_generation}. The remaining 570 samples were used exclusively for testing.
\end{itemize}

By comparing models trained on these two datasets, we aim to assess the capability of the $\phi$ML approach to extract meaningful patterns from limited data, a critical advantage when relying on expensive experimental fatigue testing.

\subsection{Training procedure}
\label{sec:training}

The machine learning models are trained using the Adam optimizer with a learning rate of \( 1 \times 10^{-4} \). The dataset is divided into training, validation, and testing sets to ensure robust model evaluation. During each training step, forward propagation is performed on mini-batches from the training set to generate predictions for cumulative fatigue life. The total custom loss \( \mathcal{L}_{\text{total}} \), integrating both data-driven and physics-based components, is then computed. Backpropagation is applied to update the model parameters using the Adam optimizer. 
The validation set is used to monitor generalization performance, ensuring appropriate hyperparameter tuning and detecting overfitting. Specifically, validation loss is tracked to assess training stability, while an early stopping criterion is implemented based on a predefined loss tolerance threshold \( \mathcal{L}_{\text{total}} \). Training is terminated if the total loss dropped below this threshold or if no significant improvement in error reduction is observed over a large number of epochs. This approach prevents unnecessary iterations and minimizes the risk of overfitting.
Throughout the training process, the total loss and its individual components are continuously monitored to assess convergence and maintain a balance between data fidelity and imposed physical constraints. This followed procedure ensures that the model not only minimizes prediction error but also adheres to the underlying physics governing concrete fatigue and loading sequence effects under non-uniform loading, enhancing predictive realism and generalizability. 
To illustrate the training duration, the machine learning model trained on the large dataset, as detailed in the subsequent sections, required 3780~seconds (63~minutes) on a standard PC equipped with 16 GB of RAM.

\section{Results and discussion on the accuracy of the $\phi$ML approach}
\label{sec:results_discussion}

This section presents a detailed evaluation of the accuracy and predictive capability of the developed machine learning models, considering different data availability scenarios and the impact of physics-based constraints.
Two distinct cases are analyzed, leveraging the datasets introduced in Sec.~\ref{sec:data_generation}. The first case explores a purely data-driven neural network trained on 70\,\% of the full dataset (630 samples), serving as a baseline to assess model performance when sufficient training data is available. This scenario represents an ideal setting where the model can learn complex fatigue lifetime patterns directly from the data.
The second case focuses on the challenge of limited data availability by training a physics-based neural network ($\phi$ML) on a small dataset of only 60 samples, as illustrated in Fig.~\ref{f:data_generation}. Unlike the first case, where extensive data allows the neural network to infer relationships purely from data patterns, this scenario reflects real-world experimental constraints, where obtaining large fatigue datasets is often impractical. To overcome the limitations of small data, the physics-based model incorporates physical constraints that guide the learning process, ensuring consistency with known fatigue behavior. For comparison, a purely data-driven neural network is also trained on the same small dataset, illustrating the advantage of integrating physical knowledge in enhancing predictive accuracy and generalization under data-scarce conditions.

\begin{figure*}[!b]
\centerline{
\includegraphics[width=1.0\textwidth]{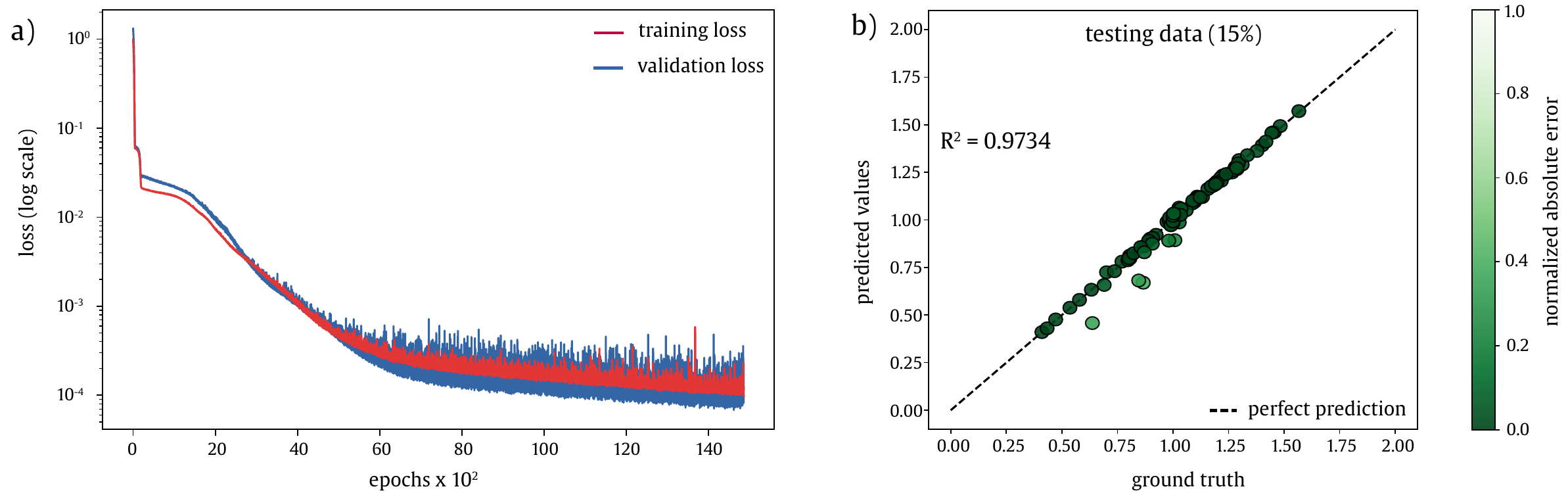}}
\caption{Results of the data-driven machine learning model trained with a large dataset (case 1): a) training history showcasing the loss function progression;\; b) prediction results for the testing dataset}
\label{f:large_data_case}
\end{figure*}

\subsection{Purely data-driven ML model trained on the large dataset}

To assess the predictive capability of the purely data-driven machine learning approach in estimating the remaining fatigue lifetime for the two different loading range scenarios (i.e. H-L and L-H), a feed-forward neural network is developed as described in Sec.~\ref{sec:NN_architecture}.
The dataset, consisting of 630 samples, is randomly split into three subsets to ensure a robust evaluation of the model’s generalization performance. Specifically, 70\,\% of the data is allocated for training, during which the model parameters are iteratively updated using back-propagation (Sec.~\ref{sec:training}). A validation set comprising 15\,\% of the data is used to monitor model performance during training, aiding in hyperparameter tuning, detecting potential overfitting, and enabling early stopping. The remaining 15\,\% of the data is reserved as a test set and remains completely unseen by the model during training, allowing for an unbiased evaluation of its final predictive capability.

The model’s training progress is illustrated in Fig.~\ref{f:large_data_case}a, hereby the evolution of the training and validation loss curves is plotted over epochs. Both curves exhibit a similar trend, with no significant divergence between them, indicating that the model successfully generalizes without severe overfitting. The training process is terminated based on the early stopping criterion, as no further significant reduction in validation loss is observed, suggesting that additional training iterations would not yield meaningful improvements.

To further analyze the model’s predictive accuracy, Fig.~\ref{f:large_data_case}b depicts the comparison between observed (i.e. ground truth) and predicted values for testing the unseen dataset. The corresponding R$^2$ score (R$^2 \in [0,1]$ score: coefficient of determination to evaluate performance) is also evaluated, providing a quantitative measure of the quality of the model fit. The results indicate that the model achieves a relatively high accuracy, with predictions closely following the observed values. 
However, a few outlier points are noticeable. These discrepancies suggest that while the model effectively captures the general trend in the data, certain regions of the input space may require additional refinement, possibly due to insufficient data density in those regions or the presence of complex, nonlinear interactions.

While this example demonstrates the effectiveness of data-driven machine learning in fatigue life prediction using a large dataset, in practical applications, especially for fatigue-related problems, data availability is often severely limited. Experimental fatigue testing is costly and time-consuming, and for complex phenomena such as the effect of loading sequence, only a small number of experimental data points are typically available. This poses a significant challenge for purely data-driven models, which rely on large datasets to generalize well. Therefore, an alternative approach is necessary to enhance predictive accuracy even with limited data. The next subsection addresses this challenge by integrating physical constraints into the machine learning model, ensuring improved generalization and reliability in data-scarce scenarios.

\begin{figure*}[!b]
\centerline{
\includegraphics[width=1.0\textwidth]{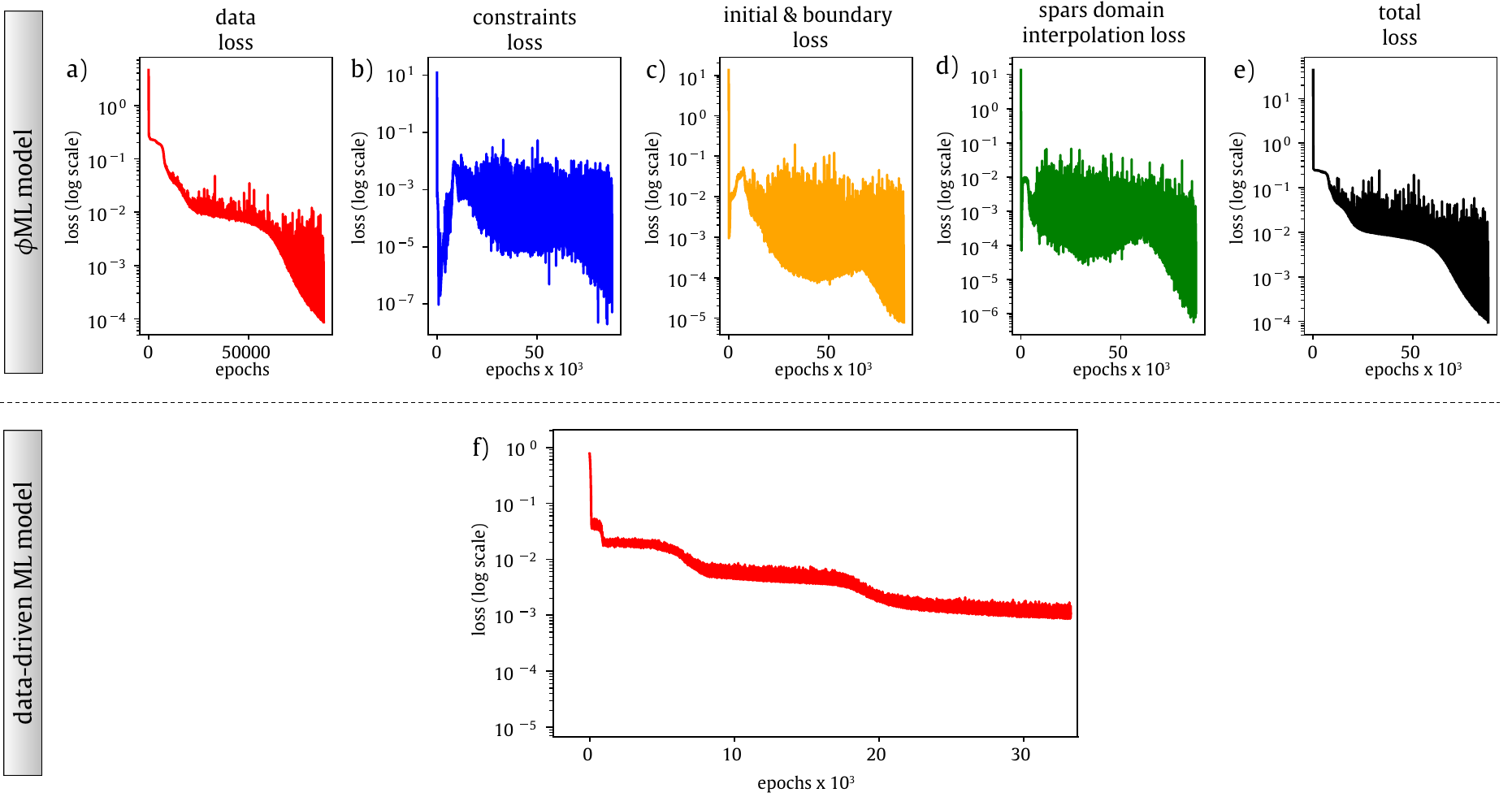}}
\caption{
Training history depicting the progression of the loss function for the data-driven ML and $\phi$ML models trained with a small dataset (case~2):
a–e) breakdown of the loss function components for the $\phi$ML model:
a) data loss;\;
b) constraints loss;\;
c) initial and boundary condition loss;\;
d) sparse domain interpolation loss;\;
e) total loss;\;
f) loss function history for the data-driven ML model
}
\label{f:loss_small_data_NN_phyNN}
\end{figure*}

\begin{figure*}[!t]
\centerline{
\includegraphics[width=1.0\textwidth]{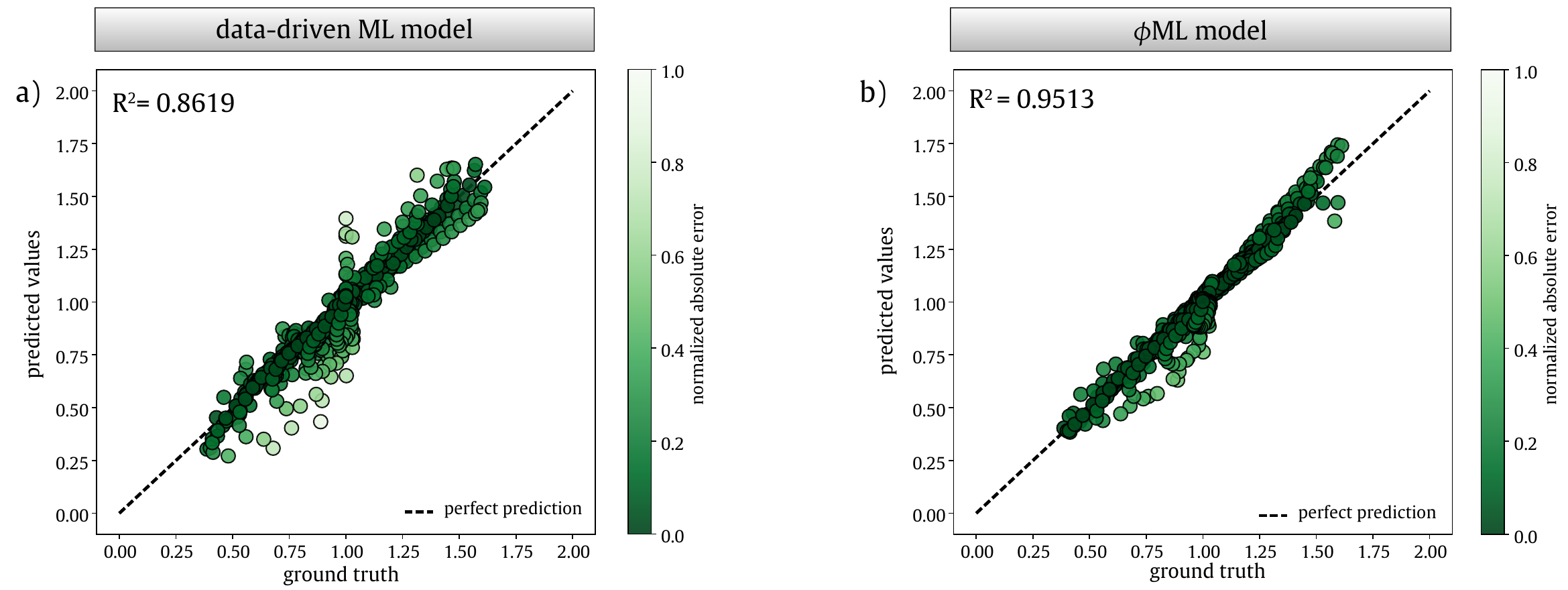}}
\caption{ 
Results of the machine learning models trained with a small dataset (case~2):
a) prediction results of the data driven ML model for the unseen testing dataset;\;
b) prediction results of the $\phi$ML model for the unseen testing dataset
}
\label{f:R_plots_small_data_NN_phyNN}
\end{figure*}

\begin{figure*}[!b]
\centerline{
\includegraphics[width=1.0\textwidth]{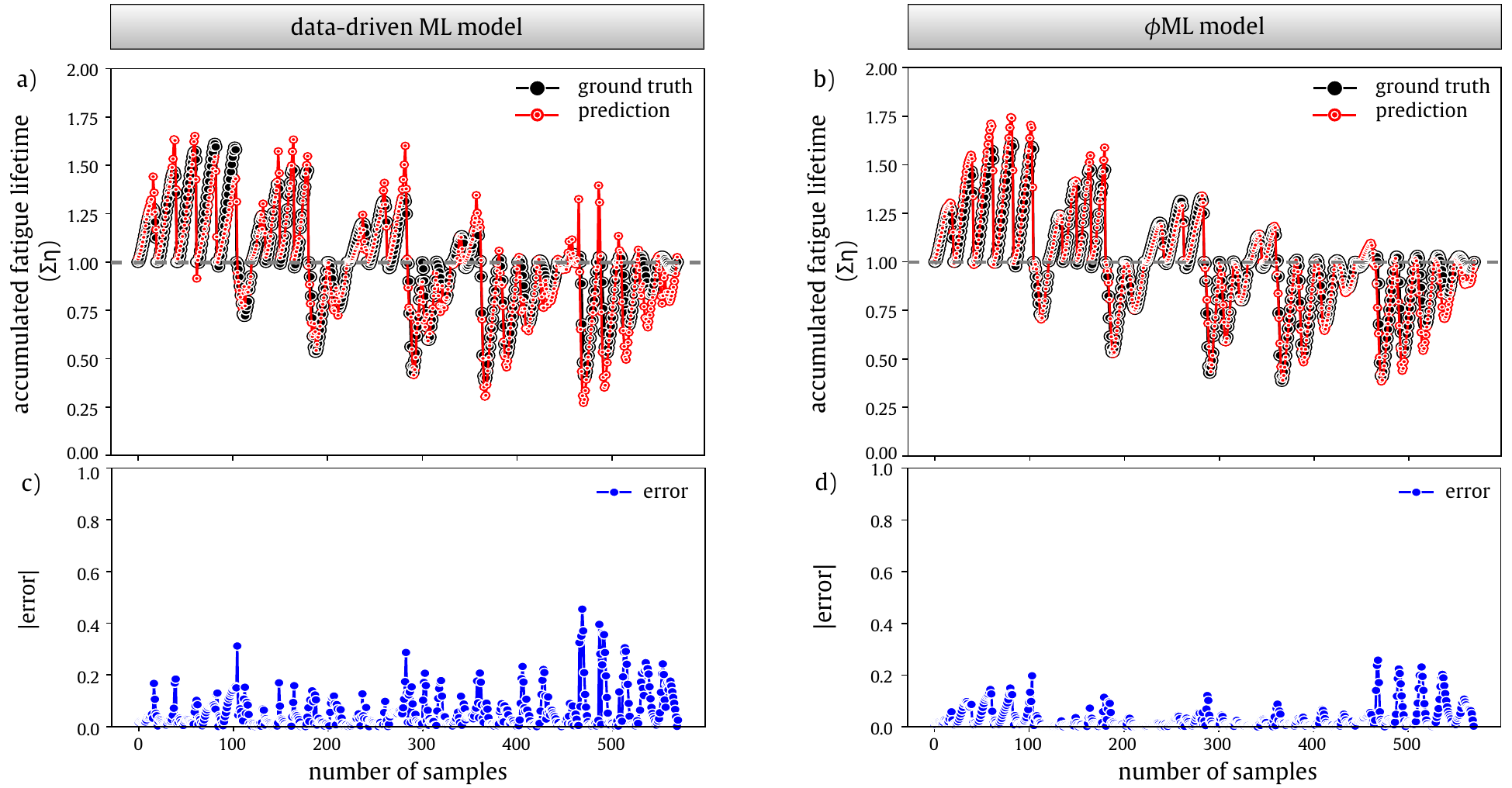}}
\caption{
Prediction results for unseen data points (samples) from the data driven ML and $\phi$ML models trained with a small dataset (case 2 and case 3, respectively):
a) comparison between the true values and predictions by the data driven ML model;\;
b) comparison between the true values and predictions by the $\phi$ML model;\;
c) absolute error values for the data driven ML model predictions;\;
d) absolute error values for the $\phi$ML model predictions
}
\label{f:error_small_data_NN_phyNN}
\end{figure*}

\subsection{Physics-based ($\phi$ML) model trained on the small dataset} \label{sec:phyNN_results_discussion}
In this case, the small dataset described in Sec.~\ref{sec:data_generation}, consisting of only 60 samples out of the total 630, is used to train the developed physics-based machine learning model. The goal here is to evaluate the capability of $\phi$ML in achieving accurate predictions and improved generalization under data-scarce conditions. The remaining 570 samples are reserved for testing to assess the accuracy of the prediction of unseen data.
For comparison, a purely data-driven FFNN is also trained using the same small dataset. This allows for a direct evaluation of the improvements gained through the introduction of physics-based constraints. The comparison between these models provides insights into the extent to which physics augmentation can mitigate the limitations of small datasets, a crucial factor in practical applications where obtaining extensive fatigue data is challenging.

Figure~\ref{f:loss_small_data_NN_phyNN} illustrates the training progress of both models through the evolution of their loss functions. For the $\phi$ML model, the first row of Fig.~\ref{f:loss_small_data_NN_phyNN} provides a detailed breakdown of its loss components, including data loss (fitting to training data), constraint loss (enforcing physical relationships), initial and boundary condition loss (ensuring meaningful extrapolation), and sparse-domain interpolation loss (improving predictions in underrepresented regions). The total loss is computed as the sum of these contributions.
The weight parameters for these loss components were tuned as follows: the constraints loss weight is set to $w_\mathrm{const} = 0.5$, the boundary condition loss weight is set to $w_\mathrm{bound} = 1.0$, and the sparse-domain interpolation loss weight is set to $w_\mathrm{spars} = 2.0$. These weights reflect the relative importance assigned to different physics-based regularizations, and strike a balance between data-driven learning and physical enforcement.

One key observation from Fig.~\ref{f:loss_small_data_NN_phyNN} is the less smooth loss evolution of the $\phi$ML model compared to the purely data-driven counterpart in Fig.~\ref{f:loss_small_data_NN_phyNN}f and the case of large-data training in Fig.~\ref{f:large_data_case}a. This behavior arises due to the additional physics-based constraints, which introduce a trade-off between minimizing the data loss and satisfying physical constraints. Unlike purely data-driven models, where the loss decreases smoothly as the model optimizes its fit to data, the physics-based model must simultaneously learn from data while satisfying multiple physical constraints, leading to occasional fluctuations in the loss evolution.
Nevertheless, despite the more complex training dynamics, the final loss values indicate that the $\phi$ML model successfully balances all loss components while maintaining stability during training.

To evaluate the predictive accuracy of the developed $\phi$ML model, the predicted values are compared with the ground-truth values for both the $\phi$ML and the purely data-driven ML models trained on a small dataset. The comparison, shown in Fig.~\ref{f:R_plots_small_data_NN_phyNN}, includes the evaluation of the R$^2$ score, which serves as a key metric for assessing model performance.
Considering the prediction of the unseen testing dataset consisting of 570 samples, a remarkable improvement in the predictive capability of the $\phi$ML model (\ref{f:R_plots_small_data_NN_phyNN}b) compared to the data-driven ML model (\ref{f:R_plots_small_data_NN_phyNN}a) is observed. These results highlight the significance of incorporating physical constraints into the machine learning framework. The $\phi$ML model, with R$^2$ score of 0.9513, demonstrates a much higher accuracy than the ML model, which achieves R$^2$ of 0.8619. This is evident in Fig.~\ref{f:R_plots_small_data_NN_phyNN}, where the predictions of $\phi$ML agree better with the actual values, while the ML model shows larger deviations.
It is noteworthy that the $\phi$ML model still holds further potential for improvement. This could be achieved by fine-tuning the hyperparameters of the custom loss function or by incorporating additional physical constraints, which would allow the model to better capture the underlying physics of the fatigue process.

To gain additional insight into the models' performance, Fig.~\ref{f:error_small_data_NN_phyNN} presents a comparison of each data point, along with the visualization of the corresponding errors. The first-row figures in this plot show the accumulated fatigue lifetime ($\sum \eta$) predictions, where the $\phi$ML model (Fig.~\ref{f:error_small_data_NN_phyNN}b) demonstrates superior performance compared to the ML model (Fig.~\ref{f:error_small_data_NN_phyNN}a). This is consistent with the error plots presented in Figs.~\ref{f:error_small_data_NN_phyNN}c and \ref{f:error_small_data_NN_phyNN}d, where the error values for the $\phi$ML are generally smaller, indicating that the $\phi$ML model is better at predicting the fatigue lifetime with less deviation from the ground truth values.
These results provide valuable insights into the trade-offs between purely data-driven and physics-based learning, offering guidance for optimizing neural network-based fatigue lifetime predictions under non-uniform loading conditions.

Additionally, the calculated average of the ground truth and predicted values of all samples is symmetrically distributed around 1.0 (highlighted as a gray dashed horizontal line), suggesting that both extensions and reductions of fatigue lifetime, observed for the L-H and H-L loading sequences, are of similar magnitudes. This observation leads to {\it an important question}: \textbf{Could the effects of the loading sequence cancel each other out in realistic loading scenarios involving a very large number of loading jumps}? This intriguing question is explored further in Section \ref{sec:generalization_validation}, where more complex loading scenarios involving multiple loading jumps are considered, which may offer a more comprehensive understanding of the interplay between loading sequences and their cumulative impact on the fatigue lifetime of materials.

%--which idealizes the uniaxial stress state into a 1D fatigue model--

\subsection{Computational efficiency}
\label{sec:comp_eff}

One of the most significant advantages of the developed machine learning models is their remarkable computational efficiency. Traditional fatigue simulations, even in simplified one-dimensional models, require cycle-by-cycle evaluations, making them computationally expensive. In the present study, the conventional cycle-by-cycle fatigue simulation requires approximately 2520~seconds (42~ minutes) to simulate $10^5$ load cycles on a standard PC with 16 GB of RAM. While this is a computationally feasible approach for simple cases, the costs for more complex simulations, especially those involving multiscale interactions or multiphysics coupled effects in the area of high-cycle fatigue up to millions of load cycles, increase dramatically, leaving it unattainable so far.
In contrast, the machine learning model achieves an instantaneous prediction of the remaining fatigue lifetime, delivering results in less than one second. As shown in Table~\ref{tab:comp_eff}, this represents a remarkable reduction in computational time without sacrificing predictive accuracy, achieving a speedup of over 2500 times. Such efficiency is particularly promising for applications requiring real-time fatigue assessment, large-scale parametric studies, or optimization tasks, where conventional simulation-based approaches would be prohibitively expensive.

\begin{table*}[!h]
\renewcommand{\arraystretch}{1.2} 
 \begin{center}
 \caption{Comparison of computational efficiency per prediction: traditional cycle-by-cycle fatigue simulation vs. $\phi$ML model}\label{tab:comp_eff}
 \begin{tabular}{ m{5.6cm} m{6.2cm} m{4.0cm} }
    \toprule
\textbf{Type of the model} & \textbf{Number of load cycles to failure} & \textbf{Computational time} \\ \hline
Traditional fatigue simulation & $10^5$ cycles & 2520 seconds (42 minutes) \\ \hline
$\phi$ML & $10^5$ cycles & \textbf{< 1 second} \\ 
  \bottomrule
 \end{tabular}
 \end{center}
\end{table*}

\section{Generalization for fatigue lifetime prediction of loading scenarios with multiple jumps}
\label{sec:generalization_validation}

\subsection{$\phi$ML based prediction of accumulated fatigue damage}
\begin{figure*}[t]
\centerline{
\includegraphics[width=1.0\textwidth]{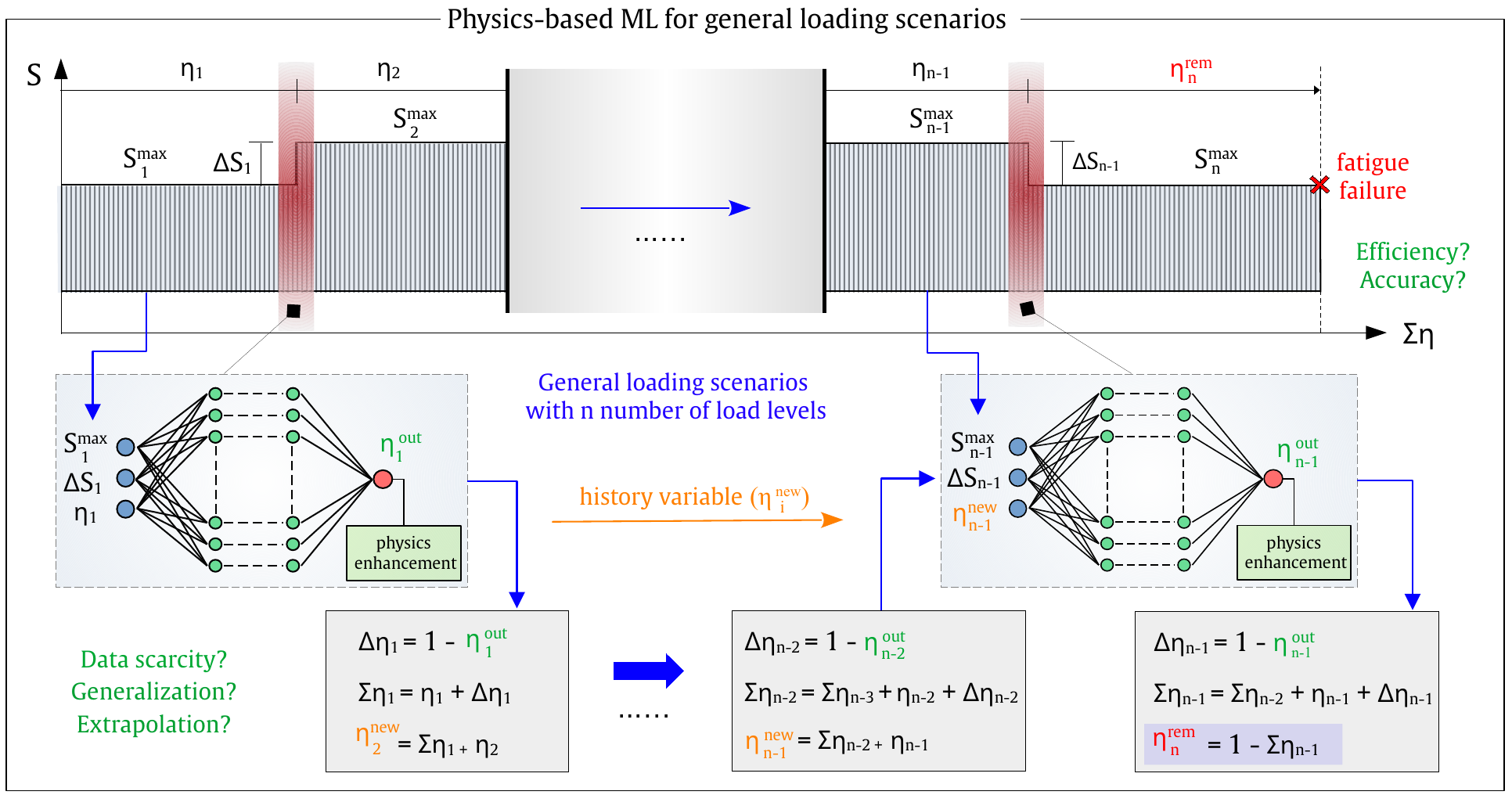}}
\caption{
Schematic representation illustrating the generalization capability of the developed $\phi$ML model for complex loading scenarios involving multiple loading jumps}
\label{f:general_concept_lifetime_prediction}
\end{figure*}

To enable accurate fatigue lifetime prediction under complex loading scenarios involving multiple load jumps, the developed $\phi$ML is integrated into an algorithm designed to estimate the remaining fatigue lifetime adaptively upon each transition between loading levels. This framework generalizes the model’s applicability to real-world conditions, where loading histories are rarely constant and often involve non-uniform sequences of cyclic stress amplitudes (see Fig.~\ref{f:general_concept_lifetime_prediction}).

Building upon the concept introduced in~\cite{BAKTHEER2019}, the proposed algorithm accounts for cumulative fatigue damage by incorporating corrections to the classical linear damage accumulation hypothesis, i.e. Palmgren-Miner (P-M) rule. This P-M rule provides a baseline estimation in which the total consumed fatigue life is the sum of individual damage fractions from different loading ranges. However, this approach neglects the interactions between consecutive load levels, which are known to significantly affect fatigue damage accumulation as described in Sec.~\ref{sec:experimental_background}. To overcome this limitation, a correction term is introduced in the current framework for each transition between the load levels that are predicted by the $\phi$ML model (see Fig.~\ref{f:general_concept_lifetime_prediction}), ensuring more realistic predictions of fatigue lifetime.
The cumulative consumed fatigue lifetime under a loading scenario consisting of $n$ varying loading levels is expressed as:
\begin{equation} \label{eq:cumulative_damage}
    \sum \eta = \sum_{i=1}^{n} \eta_\mathrm{i} +  \sum_{\mathrm{i}=1}^{n-1} \;\Delta \eta_\mathrm{i} \; ,
\end{equation}
where fatigue failure occurs when $\sum \eta = 1$, serving as the baseline failure criterion. The damage fraction for the $\mathrm{i}$-th loading level applied for the full lifetime is given by:
\begin{equation}
\eta_\mathrm{i} = \frac{N_\mathrm{i}}{N^{\mathrm{f}}_\mathrm{i}}
\end{equation}
where $N_\mathrm{i}$ is the number of applied cycles at level $\mathrm{i}$, and $N^{\mathrm{f}}_\mathrm{i}$ represents the fatigue lifetime under constant amplitude loading at that level. To extend the classical P-M rule, the correction term $\Delta \eta_{\mathrm{i}}$ is introduced to account for the effects of load history and sequence, particularly the influence of a load jump between $\mathrm{i}$ and $\mathrm{i+1}$:
\begin{equation}
\Delta \eta_\mathrm{i} = 1 - \eta^\mathrm{out}_\mathrm{i},
\end{equation}
where $\eta^\mathrm{out}_\mathrm{i}$ represents the sum of the consumed fatigue lifetime at the point of the load jump and the predicted remaining fatigue life in case the loading history continues at $S^{\text{max}}_\mathrm{i+1}$.
The $\phi$ML model developed in Sec.~\ref{sec:neural_networks} predicts this correction term based on key load transition parameters, including the maximum stress before and after the load jump ($S^{\text{max}}_\mathrm{i}$ and $S^{\text{max}}_\mathrm{i+1}$), the load amplitude change $\Delta S^{\max}_\mathrm{i}$, and the accumulated fatigue damage state $\eta^\text{cons}_\mathrm{i}$, as described in Sec.~\ref{sec:NN_architecture}.
At each load transition, the corrected accumulated fatigue lifetime, denoted as $\eta^\text{new}_\mathrm{i}$, serves as the updated history variable reflecting the current state of fatigue damage. This value is fed into the $\phi$ML model to compute the remaining fatigue life under the subsequent loading level. The iterative application of this process across multiple load jumps enables a continuous estimation of fatigue lifetime under arbitrary variable amplitude loading scenarios.

\begin{algorithm}[!t]
\caption{Remaining fatigue lifetime prediction}
\label{alg:fatigue_life}

\begin{algorithmic}[1] 
\STATE 
% \textbf{Input:} Loading scenario with $\mathrm{n}$ load levels: $\rightarrow$ LS = 
% $\begin{bmatrix}
%     S_1^{\text{max}} & S_2^{\text{max}}
%     & ... & S_\mathrm{n-1}^{\text{max}} & S_\mathrm{n}^{\text{max}} \\ 
%     \eta_1 & \eta_2 & ... & \eta_\mathrm{n-1} & \eta_\mathrm{n}: \text{output}
% \end{bmatrix}$
%
\textbf{Input:} 
Loading scenario with $\mathrm{n}$ load levels that defined by: \\\hspace{4mm} 
- Upper load level array (1  $\rightarrow$ n):  $S^{\text{max}}$ = 
$\begin{bmatrix}
    S_1^{\text{max}} & S_2^{\text{max}}
    & ... & S_\mathrm{n-1}^{\text{max}} & S_\mathrm{n}^{\text{max}}
\end{bmatrix}$\\ \vspace{1mm} \hspace{4mm} 
- Consumed lifetime array (1  $\rightarrow$ n-1): $\eta$ = 
$\begin{bmatrix}
    \eta_1 & \eta_2 & ... & \eta_\mathrm{n - 1}
\end{bmatrix}$

\vspace{1mm}

\STATE 
\textbf{Output:} Remaining fatigue life: $\rightarrow$ $\eta_\mathrm{n} = \eta^\mathrm{rem}_\mathrm{n}$

\vspace{1mm}

    \STATE Define loading levels as $\{S_\mathrm{i}^{\text{max}}, \eta_\mathrm{i} \}$ for $\text{i} = 1, \ldots, n$  with $\eta_n = 0$ (as an initial value)

    \vspace{1mm}
    
    \FOR{$\mathrm{i} = 1$ to $\mathrm{n}-1$}
    \vspace{1mm}
        \STATE \hspace{2mm} Extract the values of the current and next level: $\rightarrow$ $ S_{\text{i}}^{\text{max}}, \eta_{\text{i}}$,  $ S_{\text{i+1}}^{\text{max}}, \eta_{\text{i+1}}$

        \vspace{1mm}
        
        \STATE \hspace{2mm} Compute the current load jump: $\rightarrow$ $\Delta S^{\max}_{\text{i}} =  
        S_{\text{i}}^{\text{max}} - S_{\text{i+1}}^{\text{max}}$

        \vspace{1mm}
        
        \STATE \hspace{2mm} Compute the current consumed fatigue life: $\rightarrow$ $\eta^\text{new}_{\text{i}} = \eta_{\text{i}}$ : (\textbf{if} $\text{i} = 1$)

        % \vspace{2mm}
        
        % \IF{$i = 1$}
        %  \vspace{1mm}
        % \STATE  $\eta^\text{new}_{\text{i}} = \eta_{\text{i}}$
        %  \vspace{1mm}
        % \ENDIF

        \vspace{1mm}
        
        \STATE \hspace{2mm} Prepare the input of $\phi$ML model:  $\rightarrow$ 
        $\phi$ML$^\text{input}_\text{~i} = [ S_{\text{i}}^{\text{max}}, 
        \Delta S^{\max}_{\text{i}}, 
        \eta^\text{new}_{\text{i}}]$

        \vspace{1mm}
        
        \STATE \hspace{2mm} Call the $\phi$ML model: $\rightarrow$
        $\phi$ML$^\text{output}_\text{~i}$  = $\eta^\text{out}_\text{i}$ =
        $FFNN (\phi\text{ML}^\text{input}_\text{~i})$ 
  
        \vspace{1mm}
        
        \STATE \hspace{2mm}
        Compute the correction of the consumed fatigue life: $\rightarrow$ $\Delta \eta_{\text{i}} = 1.0 - \eta^\text{out}_\text{i}$

        \vspace{1mm}

         \STATE \hspace{2mm} Compute the current sum of consumed fatigue life: $\rightarrow$  $\sum \eta_\text{i}$
         \vspace{1mm}

         \STATE \hspace{2mm} \textbf{if} $\text{i} = 1$ \textbf{then}
         \STATE \hspace{6mm}  $\sum \eta_\text{i} = \eta_\text{i} + \Delta \eta_\text{i}$
         \STATE \hspace{2mm} \textbf{else}
         \STATE \hspace{6mm}  $\sum \eta_\text{i} = \sum \eta_\text{i-1} + \eta_\text{i} + \Delta \eta_\text{i}$
         \STATE \hspace{2mm} \textbf{end~if}

        % \IF {$i = 1$}
        %  \vspace{0mm} 
        % \STATE \hskip2mm  $\sum \eta_\text{i} = \eta_\text{i} + \Delta \eta_\text{i}$
        %  \vspace{0mm}
        % \ELSE
        %  \vspace{0mm}
        % \STATE \hskip2mm  $\sum \eta_\text{i} = \sum \eta_\text{i-1} + \eta_\text{i} + \Delta \eta_\text{i}$
        %  \vspace{0mm}
        % \ENDIF

        \vspace{1mm}

        \STATE \hspace{2mm} Compute the remaining fatigue life: $\rightarrow$
        $\eta^{\text{rem}}_\text{i+1} = 1.0 - \sum \eta_\text{i}$
        
         \vspace{1mm}
         
        \STATE \hspace{2mm} Compute the new consumed fatigue life: $\rightarrow$
        $\eta^\text{new}_{\text{i+1}} = \sum \eta_\text{i} + \eta_{\text{i+1}}$

        \vspace{1mm}

         \STATE \hspace{2mm} \textbf{if} $\eta^\text{new}_{\text{i+1}} \geq 1$
         \textbf{then}
         \STATE \hspace{6mm}  \textbf{Break}: Fatigue failure occurs within the level $\text{i+1}$
         \STATE \hspace{2mm} \textbf{end~if}
        
        % \IF{$\eta^\text{new}_{\text{i+1}} \geq 1$}
        %  \vspace{0mm}
        %     \STATE \hspace{4mm}  \textbf{Break}: Fatigue failure occurs within the level $i+1$
        %      \vspace{0mm}
        % \ENDIF
         \vspace{1mm}
    \ENDFOR

\end{algorithmic}
\end{algorithm}

The accumulated fatigue lifetime is then updated using Eq.~(\ref{eq:cumulative_damage}) and is denoted, for instance, after $n-1$ load jumps as $\sum \eta_\mathrm{n-1}$. Based on this, the remaining fatigue lifetime at the $n$-th loading stage can be estimated as:
\begin{equation}
\eta^\mathrm{rem}_n = 1 - \sum \eta_\mathrm{n-1}
\end{equation}
This formulation ensures that the accumulated damage is consistently tracked, allowing for a continuous and adaptive fatigue lifetime estimation across multiple load transitions.
The complete algorithm, incorporating the $\phi$ML at each load transition to dynamically predict fatigue lifetime, is detailed in Algorithm~\ref{alg:fatigue_life}. This approach represents a significant advancement over traditional empirical methods, offering a physics-based, data-efficient solution for fatigue lifetime assessment under complex loading conditions.

\begin{figure*}[!t]
\centerline{
\includegraphics[width=1.0\textwidth]{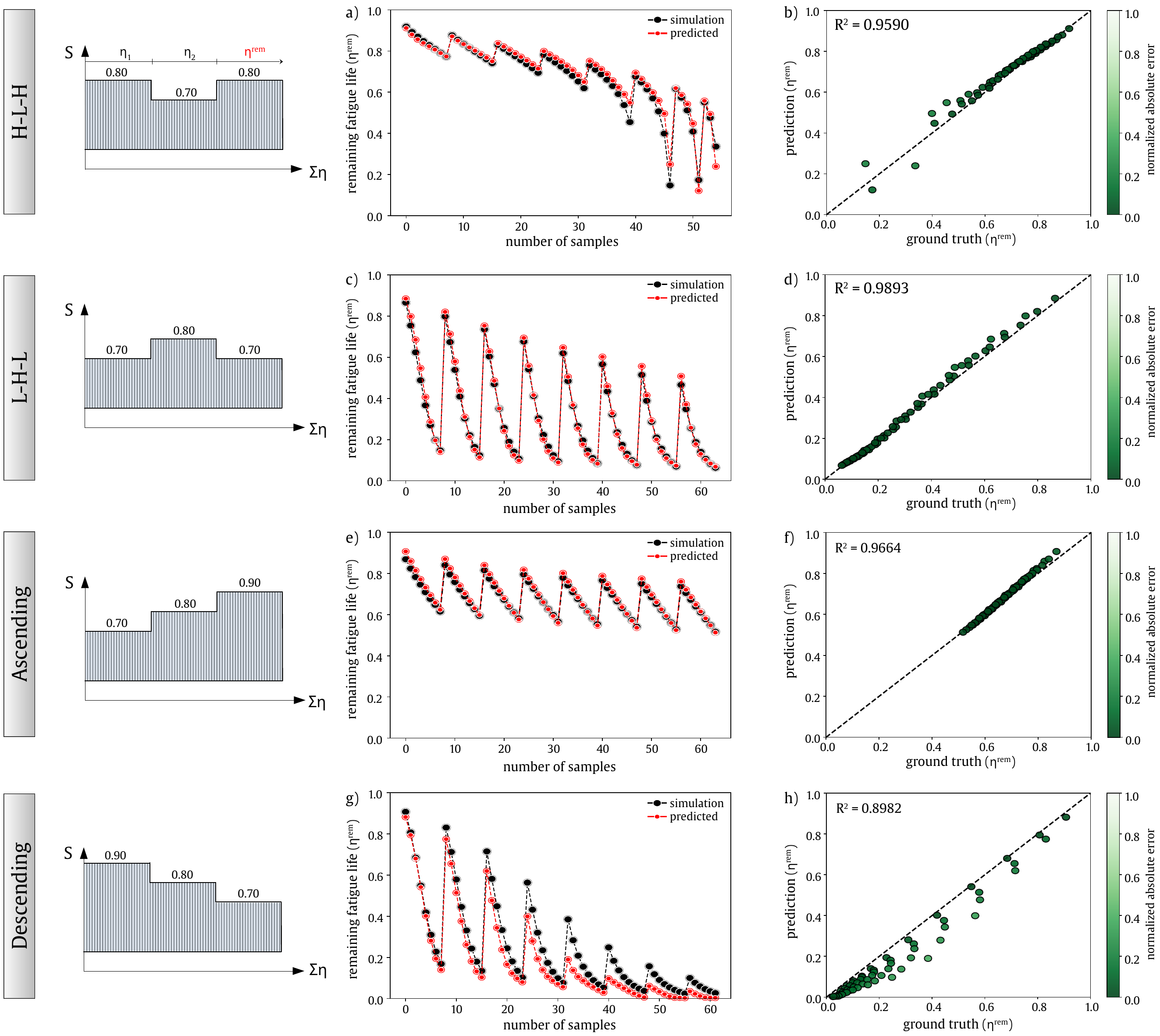}}
\caption{
Validation of the generalized $\phi$ML model for various loading scenarios with three different loading ranges in predicting the remaining fatigue lifetime:
a, c, e, g) comparison of true values and $\phi$ML model predictions for each specimen;\;
b, d, f, h) $\phi$ML model prediction results along with the corresponding R$^2$ scores
}
\label{f:three_level_validation}
\end{figure*}

\subsection{Validation utilizing three-level loading scenarios}

To assess the predictive capability of the developed $\phi$ML-based algorithm for fatigue lifetime under multiple loading levels, four different loading scenarios consisting of three load levels are examined: H-L-H, L-H-L, ascending, and descending sequences (see Fig.~\ref{f:three_level_validation}). 

Various combinations of consumed fatigue lifetime fractions $\eta_1$ and $\eta_2$ are evaluated for each scenario, corresponding to the first and second loading ranges, respectively. The values of $\eta_1$ and $\ eta_2$ systematically vary between 0.05 and 0.4 to cover a broad range of possible loading conditions. Cases where fatigue failure occurred before the second load jump are excluded from the analysis, as these effectively represent only two-level loading scenarios and do not contribute to evaluating the model’s performance for multiple loading transitions.
The remaining fatigue lifetime in the third loading range is then predicted using the developed $\phi$ML-based algorithm (Algorithm~\ref{alg:fatigue_life}) and compared against results from high-fidelity cycle-by-cycle fatigue simulations based on the anisotropic fatigue model described in Sec.~\ref{sec:fatigue_model}. The prediction accuracy is assessed through direct comparisons between the predicted and reference values, as shown in Fig.~\ref{f:three_level_validation}.

Figs.~\ref{f:three_level_validation}a, \ref{f:three_level_validation}c, \ref{f:three_level_validation}e, and \ref{f:three_level_validation}g illustrate the correlation between predicted and ground-truth values of the remaining fatigue lifetime for each loading scenario. Additionally, Figs.~\ref{f:three_level_validation}b, \ref{f:three_level_validation}d, \ref{f:three_level_validation}f, and \ref{f:three_level_validation}h visualize the prediction accuracy using R$^2$ scores, which range from 0.8982 to 0.9893, indicating reliable predictive capability across different loading sequences.
These results show that, {\it despite being trained solely on two-level loading cases}, the proposed $\phi$ML model successfully \textbf{generalizes to more complex three-level loading scenarios}. The ability to capture cumulative damage effects and adapt to varying loading sequences suggests that the physics-based approach provides an effective framework for predicting fatigue life in non-uniform loading conditions. This highlights the robustness of the developed method and its potential applicability to real-world fatigue assessment problems.

\begin{figure*}[!t]
\centerline{
\includegraphics[width=1.0\textwidth]{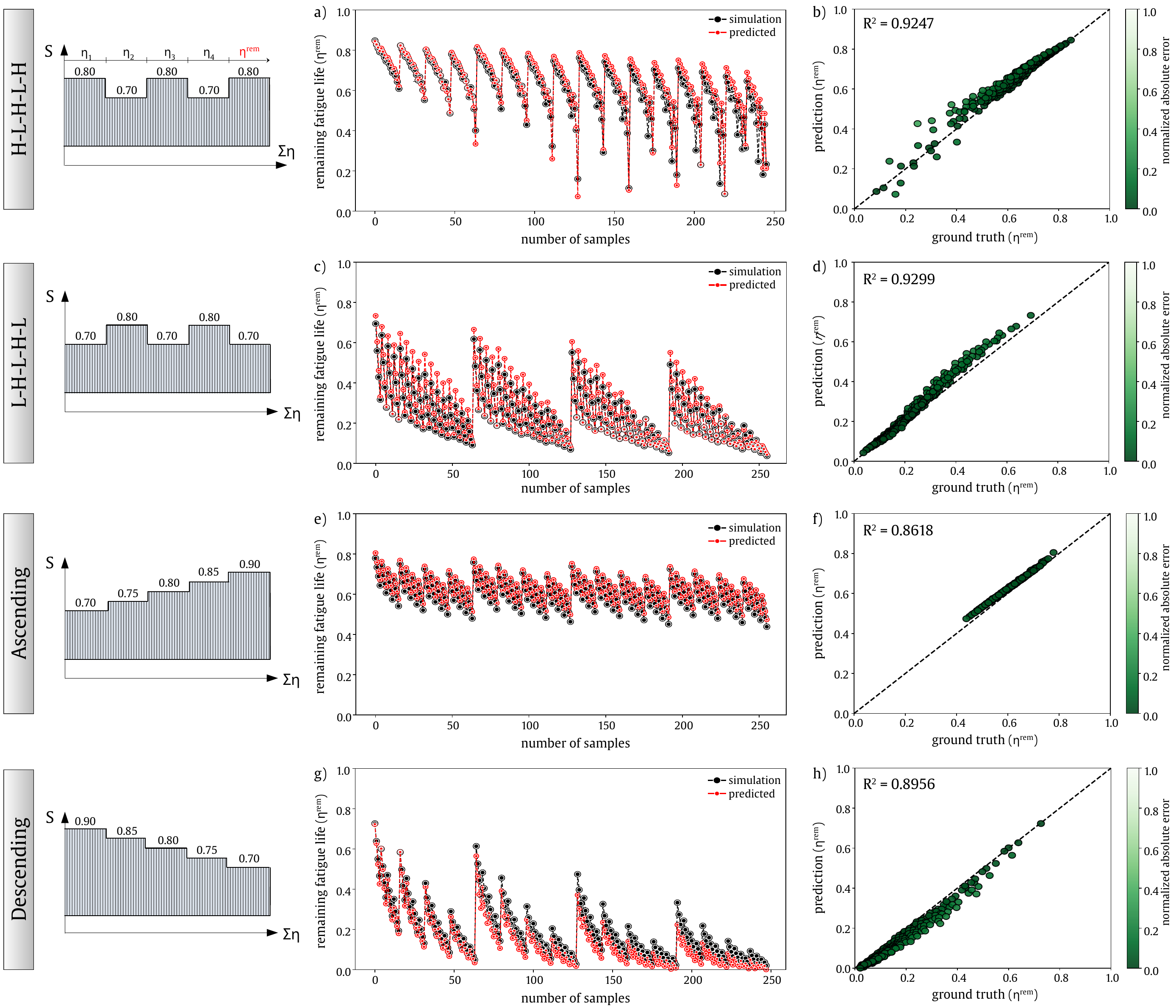}}
\caption{Validation of the generalized $\phi$ML model for various loading scenarios with five different loading ranges in predicting the remaining fatigue lifetime:
a, c, e, g) comparison of true values and $\phi$ML model predictions for each specimen;\;
b, d, f, h) $\phi$ML model prediction results along with the corresponding R$^2$ scores}
\label{f:five_level_validation}
\end{figure*}

\subsection{Validation utilizing five-level loading scenarios}

To further evaluate the generality of the developed $\phi$ML-based algorithm, an additional validation step is performed using five-level loading scenarios. The considered scenarios include H-L-H-L-H, L-H-L-H-L, as well as purely ascending and descending loading sequences as depicted in Fig.~\ref{f:five_level_validation}. These cases represent even more complex fatigue histories and allow for a deeper assessment of the model’s ability to predict fatigue damage accumulation across multiple transitions.

The results across all five-level loading scenarios show relatively good agreement between the predicted and reference fatigue lifetimes, with R$^2$ scores comparable to those obtained in the three-level validation as highlighted in Figs.~\ref{f:five_level_validation}b, \ref{f:five_level_validation}d, \ref{f:five_level_validation}f, \ref{f:five_level_validation}h. Additionally, the predicted consumed fatigue lifetime values align well with those derived from high-fidelity cycle-by-cycle fatigue simulations, further confirming the effectiveness of the presented approach.
These findings confirm the capability of the developed approach to handle complex, multi-level fatigue loading with high accuracy. The observed consistency across different loading sequences suggests that the model successfully captures the nonlinear effects of load sequence history, making it a valuable tool for fatigue lifetime prediction in engineering applications.
These used scenarios represent the first step towards complex loading histories encountered in practical fatigue applications, ensuring a comprehensive validation of the proposed approach.

\subsection{Experimental validation and analysis of accumulated fatigue damage in scenarios with multiple jumps}
\label{sec:PTST_validation}

\begin{figure*}[!t]
\centerline{
\includegraphics[width=1.0\textwidth]{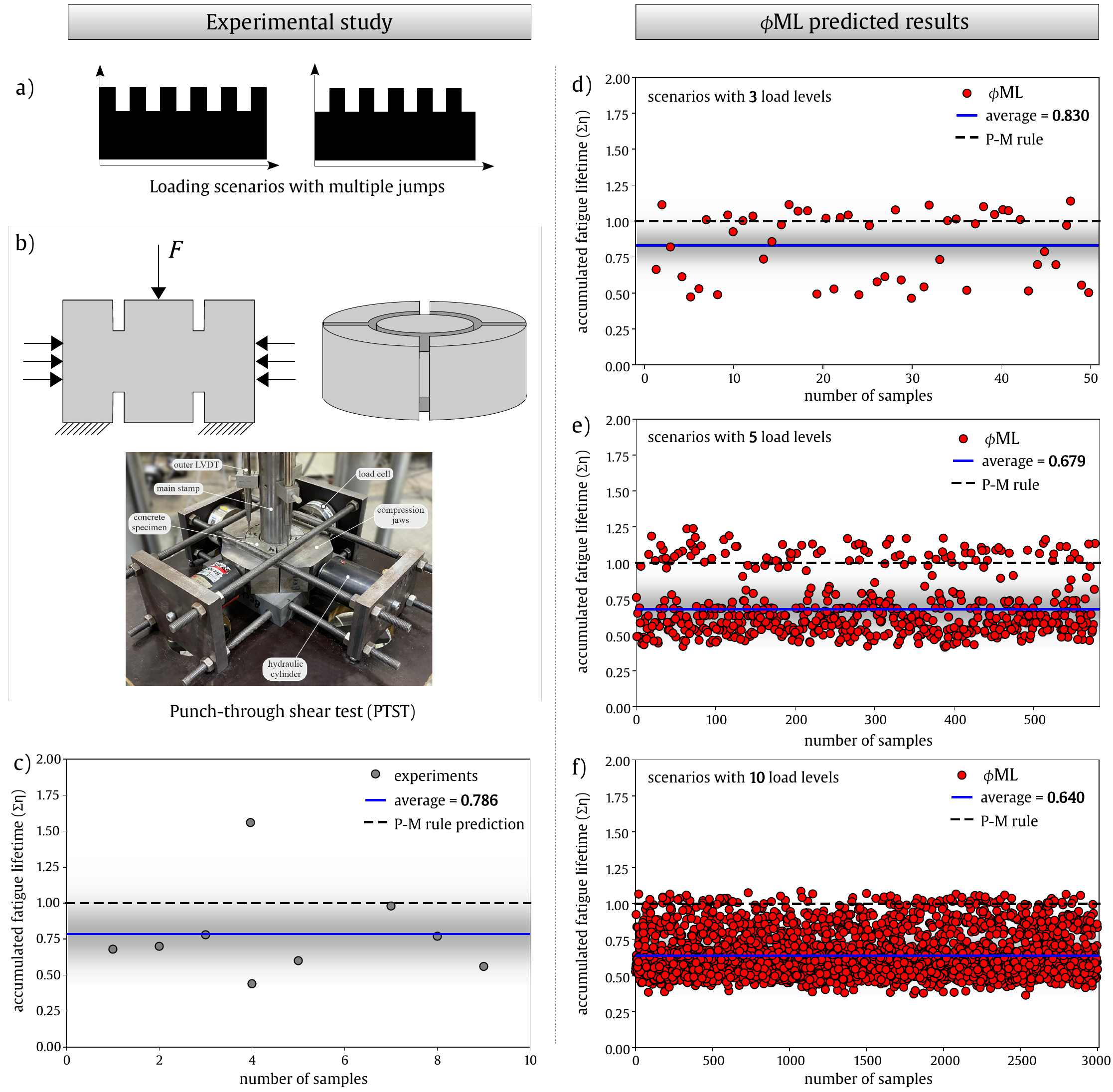}}
\caption{
Effect of multiple load jumps on the accumulated fatigue lifetime ($\sum \eta$): 
a) loading scenarios with multiple load jumps;\; 
b) punch-through shear test under combined compression-shear loading~\cite{BECKS_2024};\; 
c) experimental results for the total accumulated fatigue life;\; 
d) $\phi$ML model predictions for scenarios with 3 load levels;\; 
e) $\phi$ML model predictions for scenarios with 5 load levels;\;
f) $\phi$ML model predictions for scenarios with 10 load levels
}
\label{f:multiple_jumps_with_PTST}
\end{figure*}

Finally, to gain deeper insight into the accumulated fatigue lifetime under complex loading scenarios involving multiple jumps and to provide experimental validation for the developed $\phi$ML-based algorithm, a comparative study with recent experimental results is presented in Fig.~\ref{f:multiple_jumps_with_PTST}. This validation aims to assess the algorithm’s capability to predict fatigue lifetime in arbitrary loading scenarios and provide insights into the cumulative damage effects induced by multiple load jumps.

As highlighted in the introduction, experimental studies on the fatigue behavior of concrete under non-uniform loading conditions remain scarce in the literature. One recent study in~\cite{BECKS_2024} investigated the fatigue response of concrete under combined compression-shear loading using the punch-through shear test (PTST)~\cite{becks_mode_II, Becks_2022_Characterization, Aguilar_framcos_2023, AGUILAR_2024}, illustrated in Fig.\ref{f:multiple_jumps_with_PTST}b. This study considered multiple loading jumps with different sequences (Fig.\ref{f:multiple_jumps_with_PTST}a), systematically varying the number of load jumps between 4 and 19, with consumed fatigue lifetime fractions at each load level ranging from 0.04 to 0.2. The load jump magnitude $\Delta S^\mathrm{max}$ was consistently set to 0.2 across all tests.
The accumulated fatigue lifetimes obtained from these 9 experimental tests are depicted in Fig.~\ref{f:multiple_jumps_with_PTST}c. Despite considerable scatter in the results, a clear trend emerges, indicating that the total accumulated fatigue lifetime is generally less than the value of 1.0 (i.e., below the Palmgren-Miner (P-M) rule prediction). The average accumulated fatigue lifetime of all 9 tests was found to be $\sum \eta = 0.786$, suggesting a considerable deviation from linear damage accumulation hypothesis.

To evaluate the ability of the developed $\phi$ML-based algorithm to reproduce the observed trend of decreasing cumulative fatigue lifetime in comparison to P-M rule and to further analyze the influence of increasing load jumps, additional numerical studies are conducted. The results are shown in Figs.~\ref{f:multiple_jumps_with_PTST}d, \ref{f:multiple_jumps_with_PTST}e, and \ref{f:multiple_jumps_with_PTST}f.
Three different scenarios are investigated, each characterized by a different number of load levels (i.e., load jumps) while maintaining a constant load jump magnitude of $\Delta S^\mathrm{max} = 0.2$. The consumed fatigue life at each load range varies between 0.025 and 0.2, and the loading scenarios are generated for different combinations of these values:
\begin{itemize}
    \item Three-level loading (two load jumps): A total of 50 samples were analyzed. The predicted results using the $\phi$ML based algorithm are shown in Fig.~\ref{f:multiple_jumps_with_PTST}d. The average accumulated fatigue lifetime, indicated by the blue horizontal line, is found to be 0.830. This value is qualitatively consistent with the experimental observations in Fig.~\ref{f:multiple_jumps_with_PTST}c.

    \item Five-level loading (four load jumps): A total of 578 samples were considered. The results, presented in Fig.~\ref{f:multiple_jumps_with_PTST}e, show an average accumulated fatigue lifetime of 0.679, further confirming the trend observed experimentally.

    \item Ten-level loading (nine load jumps): To further assess the robustness of the model, an extended set of loading scenarios consisting of 3000 samples is analyzed. The results, shown in Fig.~\ref{f:multiple_jumps_with_PTST}f, reveal a similar trend with an average accumulated fatigue lifetime of 0.640.
\end{itemize}
These findings suggest a systematic decrease in the accumulated fatigue lifetime as the number of load jumps increases. However, the reduction follows a nonlinear trend, emphasizing the complexity of cumulative fatigue damage under variable amplitude loading.

The results obtained using the $\phi$ML based approach for a large number of possible loading combinations provide significant insights into the stochastic nature of fatigue loading and its impact on lifetime predictions. Notably, these findings help answer the key question raised in Sec.~\ref{sec:phyNN_results_discussion} regarding whether the lifetime extension observed in L-H scenarios and the lifetime reduction in H-L scenarios would cancel each other out with an increasing number of load jumps. The results presented in Fig.~\ref{f:multiple_jumps_with_PTST} suggest that \textbf{this is not the case}. Instead, increasing the number of load jumps leads to a systematic reduction in the accumulated fatigue lifetime, demonstrating a long-term negative effect of load history variations.

While the experimental results in Fig.~\ref{f:multiple_jumps_with_PTST}c qualitatively support this conclusion, the limited number of data points prevents definitive experimental validation. However, the developed $\phi$ML framework enables a more comprehensive exploration of these effects, as demonstrated in Figs.~\ref{f:multiple_jumps_with_PTST}d, \ref{f:multiple_jumps_with_PTST}e, and \ref{f:multiple_jumps_with_PTST}f. This suggests that further experimental investigations with a larger dataset would be valuable to solidify these conclusions.

To underscore the significant potential of physics-based machine learning $\phi$ML in fatigue analysis, consider the study depicted in Fig.~\ref{f:multiple_jumps_with_PTST}f, which encompasses 3000 distinct loading scenarios. Conducting such an extensive analysis using traditional cycle-by-cycle fatigue simulations would be \textbf{computationally too expensive}, as highlighted in Table~\ref{tab:comp_eff} for a single simulation. In contrast, the developed approach based on $\phi$ML, which only needs to be trained once (with a {\it training time of 63 minutes}), provides predictions for all \textbf{3000 scenarios in about 10~seconds}. This remarkable efficiency not only accelerates the analysis process but also enables the exploration of more complex, real-world fatigue loading conditions beyond the scope of this study.

\section{Conclusions and outlook}

The study highlighted the significant potential of the physics-based machine learning ($\phi$ML) in predicting fatigue lifetime under complex loading scenarios. Integrating the developed $\phi$ML model into a comprehensive algorithm enabled continuous estimation of both the remaining fatigue lifetime across different load levels and the total fatigue lifetime. This approach offers a promising alternative to purely data-driven methods, particularly in scenarios where data is scarce.
The algorithm has demonstrated reliable predictions of fatigue lifetime under diverse loading scenarios, including those not encountered during training. Notably, the $\phi$ML algorithm's predictions aligned with recent experimental observations, revealing a reduction in fatigue lifetime with an increased number of loading jumps. This trend was not captured by the traditional linear Palmgren-Miner (P-M) rule. This observation suggested that the P-M rule may provide overly conservative estimates of damage accumulation and fatigue life.

The efficiency and flexibility of the $\phi$ML model make it a valuable tool for integration within digital twins. Its adaptability enables the incorporation of advanced multiscale and multiphysics considerations, extending its applicability to complex fatigue phenomena in concrete and other materials. Moreover, its integration into high-fidelity multiscale fatigue modeling presents an opportunity to revolutionize fatigue lifetime prediction by significantly accelerating simulations while preserving physical consistency. This advancement paves the way for practical applications in structural health monitoring, probabilistic lifetime assessment, and real-time decision-making for fatigue-critical infrastructure and materials.
%
% In conclusion, the developed $\phi$ML-based approach represents a major advancement in fatigue lifetime prediction, combining high predictive accuracy with practical applicability. Future research should focus on further validating this approach under a broader range of loading conditions and integrating it into comprehensive structural health monitoring systems.

%\section*{CRediT authorship contribution statement}
%\textbf{Abedulgader Baktheer:} 
%Conceptualization, 
%Methodology,
%Investigation,
%Validation,
%Data Curation,
%Software,
%Writing -- Original Draft,
%Visualization.
%
%\textbf{Fadi Aldakheel:}
%Methodology, 
%Supervision,
%Visualization,
%Project administration,
%Writing -- Review \& Editing,

%\section*{Acknowledgment}
% The work was supported by the German Research Foundation (DFG) in the framework of the joint project number (XXXXXX). This support is gratefully acknowledged.
%The first author acknowledges the valuable discussions with Ahmad Delpasand from the Institute of Structural Concrete at RWTH Aachen University regarding machine learning algorithms.

\bibliography{manuscript}

\begin{thebibliography}{10}
\expandafter\ifx\csname url\endcsname\relax
  \def\url#1{\texttt{#1}}\fi
\expandafter\ifx\csname urlprefix\endcsname\relax\def\urlprefix{URL }\fi
\expandafter\ifx\csname href\endcsname\relax
  \def\href#1#2{#2} \def\path#1{#1}\fi

\bibitem{li2025arresting}
L.~Li, S.~Wang, T.~Chen, E.~Chatzi, H.~Heydarinouri, E.~Ghafoori, Arresting
  fatigue cracks in steel plates using prestressed bonded fe-sma strips:
  Analytical prediction and experimental validation, Thin-Walled Structures
  (2025) 112971.

\bibitem{Baktheer_2024_stress}
A.~Baktheer, C.~Goralski, J.~Hegger, R.~Chudoba, Stress configuration-based
  classification of current research on fatigue of reinforced and prestressed
  concrete, Structural Concrete (2024).
\newblock \href {https://doi.org/10.1002/suco.202300667}
  {\path{doi:10.1002/suco.202300667}}.

\bibitem{chen2024effect}
H.~Chen, S.~Liu, P.~Wang, X.~Wang, Z.~Liu, F.~Aldakheel, Effect of grain
  structure on fatigue crack propagation behavior of 2024 aluminum alloy under
  different stress ratios, Materials \& Design 244 (2024) 113117.

\bibitem{leuders2013mechanical}
S.~Leuders, M.~Th{\"o}ne, A.~Riemer, T.~Niendorf, T.~Tr{\"o}ster, H.~a.
  Richard, H.~Maier, On the mechanical behaviour of titanium alloy tial6v4
  manufactured by selective laser melting: Fatigue resistance and crack growth
  performance, International journal of fatigue 48 (2013) 300--307.

\bibitem{Ali_2022}
M.~Abubakar~Ali, C.~Tomann, F.~Aldakheel, M.~Mahlbacher, N.~Noii, N.~Oneschkow,
  K.-H. Drake, L.~Lohaus, P.~Wriggers, M.~Haist, Influence of moisture content
  and wet environment on the fatigue behaviour of high-strength concrete,
  Materials 15~(3) (2022).
\newblock \href {https://doi.org/10.3390/ma15031025}
  {\path{doi:10.3390/ma15031025}}.

\bibitem{van2022improving}
S.~van~den Broek, J.~Wolff, S.~Scheffler, C.~H{\"u}hne, R.~Rolfes, Improving
  the fatigue life of printed structures using stochastic variations, Progress
  in Additive Manufacturing 7~(6) (2022) 1225--1238.

\bibitem{MIARKA_2022}
P.~Miarka, S.~Seitl, V.~Bílek, H.~Cifuentes, Assessment of fatigue resistance
  of concrete: S-n curves to the paris’ law curves, Construction and Building
  Materials 341 (2022) 127811.
\newblock \href {https://doi.org/10.1016/j.conbuildmat.2022.127811}
  {\path{doi:10.1016/j.conbuildmat.2022.127811}}.

\bibitem{Seles_2021}
K.~Sele{\v{s}}, F.~Aldakheel, Z.~Tonkovi{\'{c}}, J.~Sori{\'{c}}, P.~Wriggers, A
  general phase-field model for fatigue failure in brittle and ductile solids,
  Computational Mechanics 67~(5) (2021) 1431--1452.
\newblock \href {https://doi.org/10.1007/s00466-021-01996-5}
  {\path{doi:10.1007/s00466-021-01996-5}}.

\bibitem{palmgren1924lebensdauer}
A.~Palmgren, Die {L}ebensdauer von {K}ugellagern (life length of roller
  bearings. in german), Zeitschrift des Vereines Deutscher Ingenieure (VDI
  Zeitschrift) 68 (1924) 0341--7258.

\bibitem{miner1945cumulative}
M.~Miner, et~al., Cumulative fatigue damage, Journal of applied mechanics
  12~(3) (1945) A159--A164.
\newblock \href {https://doi.org/10.1115/1.4009458}
  {\path{doi:10.1115/1.4009458}}.

\bibitem{fib_model_code_fib_2010}
Fib model code for concrete structures 2010, Document Competence Center Siegmar
  Kästl {eK}, Germany.

\bibitem{EN19922}
EN 1992-2, Eurocode 2: Design of concrete structures, Part 2: Concrete bridges
  - Design and detailing rules, European Committee for Standardisation, 2005.

\bibitem{EN19932}
EN 1993-2, Eurocode 3, Design of steel structures, Part 2: Steel Bridges,
  European Committee for Standardisation, 2006.

\bibitem{EN19931-9}
EN-1993-1-9, Eurocode 3, Design of steel structures, Part 1-9: Fatigue,
  European Committee for Standardisation, 2005.

\bibitem{marigo1985modelling}
J.~Marigo, Modelling of brittle and fatigue damage for elastic material by
  growth of microvoids, Engineering Fracture Mechanics 21~(4) (1985) 861--874.
\newblock \href {https://doi.org/10.1016/0013-7944(85)90093-1}
  {\path{doi:10.1016/0013-7944(85)90093-1}}.

\bibitem{lemaitre2005engineering}
J.~Lemaitre, R.~Desmorat, Engineering damage mechanics: ductile, creep, fatigue
  and brittle failures, Springer Science \& Business Media, 2005.

\bibitem{caggiano_2020}
A.~Caggiano, D.~S. Schicchi, S.~Yang, S.~Harenberg, V.~Malarics-Pfaff, M.~Pahn,
  F.~Dehn, E.~Koenders, A microscale approach for modelling concrete fatigue
  damage-mechanisms, in: Advances in Fracture and Damage Mechanics XVIII, Vol.
  827 of Key Engineering Materials, Trans Tech Publications Ltd, 2020, pp.
  73--78.
\newblock \href {https://doi.org/10.4028/www.scientific.net/KEM.827.73}
  {\path{doi:10.4028/www.scientific.net/KEM.827.73}}.

\bibitem{polym_6020}
E.~Martinelli, A.~Caggiano, A unified theoretical model for the monotonic and
  cyclic response of frp strips glued to concrete, Polymers 6~(2) (2014)
  370--381.
\newblock \href {https://doi.org/10.3390/polym6020370}
  {\path{doi:10.3390/polym6020370}}.

\bibitem{ULLOA_2021}
J.~Ulloa, J.~Wambacq, R.~Alessi, G.~Degrande, S.~François, Phase-field
  modeling of fatigue coupled to cyclic plasticity in an energetic formulation,
  Computer Methods in Applied Mechanics and Engineering 373 (2021) 113473.
\newblock \href {https://doi.org/10.1016/j.cma.2020.113473}
  {\path{doi:10.1016/j.cma.2020.113473}}.

\bibitem{ALESSI_2018}
R.~Alessi, S.~Vidoli, L.~{De Lorenzis}, A phenomenological approach to fatigue
  with a variational phase-field model: The one-dimensional case, Engineering
  Fracture Mechanics 190 (2018) 53--73.
\newblock \href {https://doi.org/10.1016/j.engfracmech.2017.11.036}
  {\path{doi:10.1016/j.engfracmech.2017.11.036}}.

\bibitem{Carrara_2020}
P.~Carrara, M.~Ambati, R.~Alessi, L.~{De Lorenzis}, A framework to model the
  fatigue behavior of brittle materials based on a variational phase-field
  approach, Computer Methods in Applied Mechanics and Engineering 361 (2020)
  112731.
\newblock \href {https://doi.org/10.1016/j.cma.2019.112731}
  {\path{doi:10.1016/j.cma.2019.112731}}.

\bibitem{Ghafoori_2020}
Y.~Doroudi, D.~Fernando, H.~Zhou, V.~Nguyen, E.~Ghafoori, Fatigue behavior of
  frp-to-steel bonded interface: An experimental study with a damage plasticity
  model, International Journal of Fatigue 139 (2020) 105785.
\newblock \href {https://doi.org/10.1016/j.ijfatigue.2020.105785}
  {\path{doi:10.1016/j.ijfatigue.2020.105785}}.

\bibitem{schreiber_2021}
C.~Schreiber, R.~M{\"u}ller, C.~Kuhn, Phase field simulation of fatigue crack
  propagation under complex load situations, Archive of Applied Mechanics
  91~(2) (2021) 563--577.
\newblock \href {https://doi.org/10.1007/s00419-020-01821-0}
  {\path{doi:10.1007/s00419-020-01821-0}}.

\bibitem{KHALIL_2022}
Z.~Khalil, A.~Y. Elghazouli, E.~Martínez-Pañeda, A generalised phase field
  model for fatigue crack growth in elastic–plastic solids with an efficient
  monolithic solver, Computer Methods in Applied Mechanics and Engineering 388
  (2022) 114286.
\newblock \href {https://doi.org/10.1016/j.cma.2021.114286}
  {\path{doi:10.1016/j.cma.2021.114286}}.

\bibitem{SCHRODER_2022}
J.~Schröder, M.~Pise, D.~Brands, G.~Gebuhr, S.~Anders, Phase-field modeling of
  fracture in high performance concrete during low-cycle fatigue: Numerical
  calibration and experimental validation, Computer Methods in Applied
  Mechanics and Engineering 398 (2022) 115181.
\newblock \href {https://doi.org/10.1016/j.cma.2022.115181}
  {\path{doi:10.1016/j.cma.2022.115181}}.

\bibitem{kobler2021computational}
J.~K{\"o}bler, N.~Magino, H.~Andr{\"a}, F.~Welschinger, R.~M{\"u}ller,
  M.~Schneider, A computational multi-scale model for the stiffness degradation
  of short-fiber reinforced plastics subjected to fatigue loading, Computer
  Methods in Applied Mechanics and Engineering 373 (2021) 113522.

\bibitem{lee2021review}
H.~W. Lee, C.~Basaran, A review of damage, void evolution, and fatigue life
  prediction models, Metals 11~(4) (2021) 609.

\bibitem{kirane2015microplane}
K.~Kirane, Z.~P. Ba{\v{z}}ant, Microplane damage model for fatigue of
  quasibrittle materials: Sub-critical crack growth, lifetime and residual
  strength, International Journal of Fatigue 70 (2015) 93--105.
\newblock \href {https://doi.org/10.1016/j.ijfatigue.2014.08.012}
  {\path{doi:10.1016/j.ijfatigue.2014.08.012}}.

\bibitem{BAKTHEER_2021_MS1}
A.~Baktheer, M.~Aguilar, R.~Chudoba, Microplane fatigue model {MS}1 for plain
  concrete under compression with damage evolution driven by cumulative
  inelastic shear strain, International Journal of Plasticity 143 (2021).
\newblock \href {https://doi.org/10.1016/j.ijplas.2021.102950}
  {\path{doi:10.1016/j.ijplas.2021.102950}}.

\bibitem{Baktheer_2024_FFEMS}
A.~Baktheer, S.~Esfandiari, M.~Aguilar, H.~Becks, M.~Classen, R.~Chudoba,
  Fatigue-induced stress redistribution in prestressed concrete beams modeled
  using the constitutive hypothesis of inter-aggregate degradation, Fatigue \&
  Fracture of Engineering Materials \& Structures 47~(10) (2024) 3673--3692.
\newblock \href {https://doi.org/10.1111/ffe.14388}
  {\path{doi:10.1111/ffe.14388}}.

\bibitem{hessman2023micromechanical}
P.~A. Hessman, F.~Welschinger, K.~Hornberger, T.~B{\"o}hlke, A micromechanical
  cyclic damage model for high cycle fatigue failure of short fiber reinforced
  composites, Composites Part B: Engineering 264 (2023) 110855.

\bibitem{Ueda2019}
N.~Ueda, M.~Konishi, H.~Ono, Quasi-visco-elasto-plastic constitutive model of
  concrete for fatigue simulation, FraMCoS-X Conference in Bayonne, France,
  2019.
\newblock \href {https://doi.org/10.21012/FC10.233513}
  {\path{doi:10.21012/FC10.233513}}.

\bibitem{LE_2023}
V.~T. Le, K.~M. Tran, J.~Kodikara, D.~Bodin, J.~Grenfell, H.~H. Bui, A
  two-surface contact model for dem and its application to model fatigue crack
  growth in cemented materials, International Journal of Plasticity 166 (2023)
  103650.
\newblock \href {https://doi.org/10.1016/j.ijplas.2023.103650}
  {\path{doi:10.1016/j.ijplas.2023.103650}}.

\bibitem{shojai2024micro}
S.~Shojai, P.~Schaumann, E.~Ghafoori, Micro-support effect consideration in
  fatigue analysis of corroded steel based on real surface geometry, Journal of
  Constructional Steel Research 212 (2024) 108259.

\bibitem{zhan2021machine}
Z.~Zhan, H.~Li, Machine learning based fatigue life prediction with effects of
  additive manufacturing process parameters for printed ss 316l, Int. J.
  Fatigue 142 (2021) 105941.
\newblock \href {https://doi.org/10.1016/j.ijfatigue.2020.105941}
  {\path{doi:10.1016/j.ijfatigue.2020.105941}}.

\bibitem{Son_2022}
J.~Son, S.~Yang, A new approach to machine learning model development for
  prediction of concrete fatigue life under uniaxial compression, Applied
  Sciences 12~(19) (2022).
\newblock \href {https://doi.org/10.3390/app12199766}
  {\path{doi:10.3390/app12199766}}.

\bibitem{Zhang_2022}
W.~Zhang, D.~Lee, J.~Lee, C.~Lee, Residual strength of concrete subjected to
  fatigue based on machine learning technique, Structural Concrete 23~(4)
  (2022) 2274--2287.
\newblock \href {https://doi.org/10.1002/suco.202100082}
  {\path{doi:10.1002/suco.202100082}}.

\bibitem{zhan2021data}
Z.~Zhan, W.~Hu, Q.~Meng, Data-driven fatigue life prediction in additive
  manufactured titanium alloy: A damage mechanics based machine learning
  framework, Engineering Fracture Mechanics 252 (2021) 107850.

\bibitem{RIYAR_2023}
R.~L. Riyar, Mansi, S.~Bhowmik, Fatigue behaviour of plain and reinforced
  concrete: A systematic review, Theoretical and Applied Fracture Mechanics 125
  (2023) 103867.
\newblock \href {https://doi.org/10.1016/j.tafmec.2023.103867}
  {\path{doi:10.1016/j.tafmec.2023.103867}}.

\bibitem{hamada2025advancing}
A.~Hamada, S.~Elyamny, W.~Abd-Elaziem, S.~Elkatatny, M.~A. Darwish, T.~A.
  Sebaey, A.~J{\"a}rvenp{\"a}{\"a}, K.~Vineesh, A.~H. Elsheikh, Advancing
  fatigue life prediction with machine learning: A review, Materials Today
  Communications (2025) 111525.

\bibitem{ZHANG_2022_ML_fatigue_1}
J.~Zhang, J.~Zhu, W.~Guo, W.~Guo, A machine learning-based approach to predict
  the fatigue life of three-dimensional cracked specimens, International
  Journal of Fatigue 159 (2022) 106808.
\newblock \href {https://doi.org/10.1016/j.ijfatigue.2022.106808}
  {\path{doi:10.1016/j.ijfatigue.2022.106808}}.

\bibitem{YANG_2022_ML_fatigue_2}
J.~Yang, G.~Kang, Q.~Kan, A novel deep learning approach of multiaxial fatigue
  life-prediction with a self-attention mechanism characterizing the effects of
  loading history and varying temperature, International Journal of Fatigue 162
  (2022) 106851.
\newblock \href {https://doi.org/10.1016/j.ijfatigue.2022.106851}
  {\path{doi:10.1016/j.ijfatigue.2022.106851}}.

\bibitem{app15052589}
Y.~Heider, F.~Aldakheel, W.~Ehlers, A multiscale cnn-based intrinsic
  permeability prediction in deformable porous media, Applied Sciences 15~(5)
  (2025).
\newblock \href {https://doi.org/10.3390/app15052589}
  {\path{doi:10.3390/app15052589}}.

\bibitem{montans2023}
F.~J. Montáns, E.~Cueto, K.-J. Bathe, Machine {Learning} in {Computer} {Aided}
  {Engineering}, in: T.~Rabczuk, K.-J. Bathe (Eds.), Machine {Learning} in
  {Modeling} and {Simulation}: {Methods} and {Applications}, Computational
  {Methods} in {Engineering} \& the {Sciences}, Springer International
  Publishing, Cham, 2023, pp. 1--83.
\newblock \href {https://doi.org/10.1007/978-3-031-36644-4_1}
  {\path{doi:10.1007/978-3-031-36644-4_1}}.

\bibitem{karniadakis2021piml}
G.~E. Karniadakis, I.~G. Kevrekidis, L.~Lu, P.~Perdikaris, S.~Wang, L.~Yang,
  Physics-informed machine learning 3~(6)  422--440.
\newblock \href {https://doi.org/10.1038/s42254-021-00314-5}
  {\path{doi:10.1038/s42254-021-00314-5}}.

\bibitem{ASAD_2023}
F.~As’ad, C.~Farhat, A mechanics-informed deep learning framework for
  data-driven nonlinear viscoelasticity, Computer Methods in Applied Mechanics
  and Engineering 417 (2023) 116463.
\newblock \href {https://doi.org/10.1016/j.cma.2023.116463}
  {\path{doi:10.1016/j.cma.2023.116463}}.

\bibitem{michopoulos2024scientific}
J.~G. Michopoulos, A.~Bhaduri, F.~Chinesta, E.~Cueto, D.~Liu, S.~K. Ravi, J.-X.
  Wang, Scientific machine learning for manufacturing processes and material
  systems, Journal of Computing and Information Science in Engineering 24~(11)
  (2024).

\bibitem{Cueto_2023}
E.~Cueto, F.~Chinesta, Thermodynamics of learning physical phenomena, Archives
  of Computational Methods in Engineering 30~(8) (2023) 4653--4666.
\newblock \href {https://doi.org/10.1007/s11831-023-09954-5}
  {\path{doi:10.1007/s11831-023-09954-5}}.

\bibitem{ELFALLAKIIDRISSI_2024}
M.~{El Fallaki Idrissi}, F.~Praud, F.~Meraghni, F.~Chinesta, G.~Chatzigeorgiou,
  Multiscale thermodynamics-informed neural networks (mutinn) towards fast and
  frugal inelastic computation of woven composite structures, Journal of the
  Mechanics and Physics of Solids 186 (2024) 105604.
\newblock \href {https://doi.org/10.1016/j.jmps.2024.105604}
  {\path{doi:10.1016/j.jmps.2024.105604}}.

\bibitem{di2024physics}
D.~Di~Lorenzo, V.~Champaney, C.~Ghnatios, E.~Cueto, F.~Chinesta,
  Physics-informed and graph neural networks for enhanced inverse analysis,
  Engineering Computations (2024).
\newblock \href {https://doi.org/10.1108/EC-12-2023-0958}
  {\path{doi:10.1108/EC-12-2023-0958}}.

\bibitem{manav2024phase}
M.~Manav, R.~Molinaro, S.~Mishra, L.~De~Lorenzis, Phase-field modeling of
  fracture with physics-informed deep learning, Computer Methods in Applied
  Mechanics and Engineering 429 (2024) 117104.

\bibitem{THIERCELIN_2024}
L.~Thiercelin, L.~Peltier, F.~Meraghni, Physics-informed machine learning
  prediction of the martensitic transformation temperature for the design of
  “niti-like” high entropy shape memory alloys, Computational Materials
  Science 231 (2024) 112578.
\newblock \href {https://doi.org/10.1016/j.commatsci.2023.112578}
  {\path{doi:10.1016/j.commatsci.2023.112578}}.

\bibitem{LEON_2025}
O.~León, V.~Rivera, A.~Vázquez-Patiño, J.~Ulloa, E.~Samaniego, Exploring
  energy minimization to model strain localization as a strong discontinuity
  using physics informed neural networks, Computer Methods in Applied Mechanics
  and Engineering 436 (2025) 117724.
\newblock \href {https://doi.org/10.1016/j.cma.2024.117724}
  {\path{doi:10.1016/j.cma.2024.117724}}.

\bibitem{Athanasiou_2025}
R.~Yi, D.~Georgiou, X.~Liu, C.~E. Athanasiou, Mechanics-informed, model-free
  symbolic regression framework for solving fracture problems, Journal of the
  Mechanics and Physics of Solids 194 (2025) 105916.
\newblock \href {https://doi.org/10.1016/j.jmps.2024.105916}
  {\path{doi:10.1016/j.jmps.2024.105916}}.

\bibitem{SALVATI_2022_ML_fatigue_3}
E.~Salvati, A.~Tognan, L.~Laurenti, M.~Pelegatti, F.~{De Bona}, A defect-based
  physics-informed machine learning framework for fatigue finite life
  prediction in additive manufacturing, Materials \& Design 222 (2022) 111089.
\newblock \href {https://doi.org/10.1016/j.matdes.2022.111089}
  {\path{doi:10.1016/j.matdes.2022.111089}}.

\bibitem{HALAMKA_2023_ML_fatigue_6}
J.~Halamka, M.~Bartošák, M.~Španiel, Using hybrid physics-informed neural
  networks to predict lifetime under multiaxial fatigue loading, Engineering
  Fracture Mechanics 289 (2023) 109351.
\newblock \href {https://doi.org/10.1016/j.engfracmech.2023.109351}
  {\path{doi:10.1016/j.engfracmech.2023.109351}}.

\bibitem{HE_2023__ML_fatigue_8}
G.~He, Y.~Zhao, C.~Yan, Multiaxial fatigue life prediction using
  physics-informed neural networks with sensitive features, Engineering
  Fracture Mechanics 289 (2023) 109456.
\newblock \href {https://doi.org/10.1016/j.engfracmech.2023.109456}
  {\path{doi:10.1016/j.engfracmech.2023.109456}}.

\bibitem{WANG_2023_ML_fatigue_7}
L.~Wang, S.-P. Zhu, C.~Luo, D.~Liao, Q.~Wang, Physics-guided machine learning
  frameworks for fatigue life prediction of am materials, International Journal
  of Fatigue 172 (2023) 107658.
\newblock \href {https://doi.org/10.1016/j.ijfatigue.2023.107658}
  {\path{doi:10.1016/j.ijfatigue.2023.107658}}.

\bibitem{JIANG_2024_ML_fatigue_5}
L.~Jiang, Y.~Hu, Y.~Liu, X.~Zhang, G.~Kang, Q.~Kan, Physics-informed machine
  learning for low-cycle fatigue life prediction of 316 stainless steels,
  International Journal of Fatigue 182 (2024) 108187.
\newblock \href {https://doi.org/10.1016/j.ijfatigue.2024.108187}
  {\path{doi:10.1016/j.ijfatigue.2024.108187}}.

\bibitem{gonzalez2020}
D.~Gonzalez, A.~Garcia-Gonzalez, F.~Chinesta, E.~Cueto, A data-driven learning
  method for constitutive modeling: Application to vascular hyperelastic soft
  tissues, Materials 13~(10), https://www.mdpi.com/1996-1944/13/10/2319 (2020).
\newblock \href {https://doi.org/10.3390/ma13102319}
  {\path{doi:10.3390/ma13102319}}.

\bibitem{klein2022}
D.~K. Klein, M.~Fernandez, R.~J. Martin, P.~Neff, O.~Weeger, Polyconvex
  anisotropic hyperelasticity with neural networks, J. Mech. Phys. Solids 159
  (2022) 104703.
\newblock \href {https://doi.org/10.1016/j.jmps.2021.104703}
  {\path{doi:10.1016/j.jmps.2021.104703}}.

\bibitem{tac2023a}
V.~Tac, M.~K. Rausch, F.~Sahli~Costabal, A.~B. Tepole, Data-driven anisotropic
  finite viscoelasticity using neural ordinary differential equations, Comput.
  Methods Appl. Mech. Eng. 411 (2023) 116046.
\newblock \href {https://doi.org/10.1016/j.cma.2023.116046}
  {\path{doi:10.1016/j.cma.2023.116046}}.

\bibitem{rosenkranz2024a}
M.~Rosenkranz, K.~A. Kalina, J.~Brummund, W.~Sun, M.~K\"astner, Viscoelasticty
  with physics-augmented neural networks: Model formulation and training
  methods without prescribed internal variables, Computational
  MechanicsHttps://link.springer.com/10.1007/s00466-024-02477-1 (2024).
\newblock \href {https://doi.org/10.1007/s00466-024-02477-1}
  {\path{doi:10.1007/s00466-024-02477-1}}.

\bibitem{vlassis2021}
N.~N. Vlassis, W.~Sun, Sobolev training of thermodynamic-informed neural
  networks for interpretable elasto-plasticity models with level set hardening,
  Comput. Methods Appl. Mech. Eng. 377 (2021) 113695.
\newblock \href {https://doi.org/10.1016/j.cma.2021.113695}
  {\path{doi:10.1016/j.cma.2021.113695}}.

\bibitem{masi2023}
F.~Masi, I.~Stefanou, Evolution {TANN} and the identification of internal
  variables and evolution equations in solid mechanics, Journal of the
  Mechanics and Physics of Solids 174 (2023) 105245.
\newblock \href {https://doi.org/10.1016/j.jmps.2023.105245}
  {\path{doi:10.1016/j.jmps.2023.105245}}.

\bibitem{fuhg2023}
J.~N. Fuhg, C.~M. Hamel, K.~Johnson, R.~Jones, N.~Bouklas, Modular machine
  learning-based elastoplasticity: {{Generalization}} in the context of limited
  data, Comput. Methods Appl. Mech. Eng. 407 (2023) 115930.

\bibitem{Aldakheel_2025}
F.~Aldakheel, E.~S. Elsayed, Y.~Heider, O.~Weeger, Physics-based machine
  learning for computational fracture mechanics, arXiv preprint
  arXiv:2502.09025 (2025).
\newblock \href {https://doi.org/10.48550/arXiv.2502.09025}
  {\path{doi:10.48550/arXiv.2502.09025}}.

\bibitem{Agerskov}
H.~Agerskov, N.~T. Pedersen, Fatigue life of offshore steel structures under
  stochastic loading, Journal of Structural Engineering 118~(8) (1992)
  2101--2117.
\newblock \href {https://doi.org/10.1061/(ASCE)0733-9445(1992)118:8(2101)}
  {\path{doi:10.1061/(ASCE)0733-9445(1992)118:8(2101)}}.

\bibitem{Tepfers_1977}
R.~Tepfers, C.~Frid\'{e}n, L.~Georgsson, A study of the applicability to the
  fatigue of concrete of the palmgren-miner partial damage hypothesis, Magazine
  of Concrete Research 29~(100) (1977) 123--130.
\newblock \href {https://doi.org/10.1680/macr.1977.29.100.123}
  {\path{doi:10.1680/macr.1977.29.100.123}}.

\bibitem{holmen1982fatigue}
J.~O. Holmen, Fatigue of concrete by constant and variable amplitude loading,
  ACI Journal, Special Publication 75 (1982) 71--110.
\newblock \href {https://doi.org/10.14359/6402} {\path{doi:10.14359/6402}}.

\bibitem{petkovic1990fatigue}
G.~Petkovic, R.~Lenschow, H.~Stemland, S.~Rosseland, Fatigue of high-strength
  concrete, Special Publication 121 (1990) 505--526.
\newblock \href {https://doi.org/10.14359/3740} {\path{doi:10.14359/3740}}.

\bibitem{Baktheer_2021_3}
A.~Baktheer, R.~Chudoba, Experimental and theoretical evidence for the load
  sequence effect in the compressive fatigue behavior of concrete, Materials
  and Structures 54~(2) (2021) 82.
\newblock \href {https://doi.org/10.1617/s11527-021-01667-0}
  {\path{doi:10.1617/s11527-021-01667-0}}.

\bibitem{hoff1984testing}
A.~Hoff, Testing of high strength lightweight aggregate concrete elements,
  Nordic concrete research~(3) (1984) 63--91.

\bibitem{hilsdorf1966fatigue}
H.~K. Hilsdorf, et~al., Fatigue strength of concrete under varying flexural
  stresses, in: Journal Proceedings, Vol.~63, 1966, pp. 1059--1076.
\newblock \href {https://doi.org/10.14359/7662} {\path{doi:10.14359/7662}}.

\bibitem{BAKTHEER_2021_4}
A.~Baktheer, H.~Becks, Fracture mechanics based interpretation of the load
  sequence effect in the flexural fatigue behavior of concrete using digital
  image correlation, Construction and Building Materials 307 (2021).
\newblock \href {https://doi.org/10.1016/j.conbuildmat.2021.124817}
  {\path{doi:10.1016/j.conbuildmat.2021.124817}}.

\bibitem{baktheer2018modeling}
A.~Baktheer, R.~Chudoba, Modeling of bond fatigue in reinforced concrete based
  on cumulative measure of slip, in: Computational Modelling of Concrete
  Structures, CRC Press, 2018, pp. 767--776.

\bibitem{HALM1998439}
D.~Halm, A.~Dragon, An anisotropic model of damage and frictional sliding for
  brittle materials, European Journal of Mechanics - A/Solids 17~(3) (1998) 439
  -- 460.
\newblock \href {https://doi.org/10.1016/S0997-7538(98)80054-5}
  {\path{doi:10.1016/S0997-7538(98)80054-5}}.

\bibitem{DRAGON2000331}
A.~Dragon, D.~Halm, T.~Désoyer, Anisotropic damage in quasi-brittle solids:
  modelling, computational issues and applications, Computer Methods in Applied
  Mechanics and Engineering 183~(3) (2000) 331 -- 352.
\newblock \href {https://doi.org/10.1016/S0045-7825(99)00225-X}
  {\path{doi:10.1016/S0045-7825(99)00225-X}}.

\bibitem{alliche2004damage}
A.~Alliche, Damage model for fatigue loading of concrete, International Journal
  of Fatigue 26~(9) (2004) 915--921.
\newblock \href {https://doi.org/10.1016/j.ijfatigue.2004.02.006}
  {\path{doi:10.1016/j.ijfatigue.2004.02.006}}.

\bibitem{BAKTHEER2019}
A.~Baktheer, J.~Hegger, R.~Chudoba, Enhanced assessment rule for concrete
  fatigue under compression considering the nonlinear effect of loading
  sequence, International Journal of Fatigue 126 (2019) 130 -- 142.
\newblock \href {https://doi.org/10.1016/j.ijfatigue.2019.04.027}
  {\path{doi:10.1016/j.ijfatigue.2019.04.027}}.

\bibitem{mazars1989continuum}
J.~Mazars, G.~Pijaudier-Cabot, Continuum damage theory—application to
  concrete, Journal of Engineering Mechanics 115~(2) (1989) 345--365.
\newblock \href {https://doi.org/10.1061/(ASCE)0733-9399(1989)115:2(345)}
  {\path{doi:10.1061/(ASCE)0733-9399(1989)115:2(345)}}.

\bibitem{Baktheer_2024_PFM}
A.~Baktheer, E.~Martínez-Pañeda, F.~Aldakheel, Phase field cohesive zone
  modeling for fatigue crack propagation in quasi-brittle materials, Computer
  Methods in Applied Mechanics and Engineering 422 (2024) 116834.
\newblock \href {https://doi.org/10.1016/j.cma.2024.116834}
  {\path{doi:10.1016/j.cma.2024.116834}}.

\bibitem{Oneschkow_Loehnert_2022}
N.~Oneschkow, T.~Timmermann, S.~Löhnert, Compressive fatigue behaviour of
  high-strength concrete and mortar: Experimental investigations and
  computational modelling, Materials 15~(1) (2022).
\newblock \href {https://doi.org/10.3390/ma15010319}
  {\path{doi:10.3390/ma15010319}}.

\bibitem{Kristensen2023}
P.~K. Kristensen, A.~Golahmar, E.~Martínez-Pañeda, C.~F. Niordson,
  Accelerated high-cycle phase field fatigue predictions, {European Journal of
  Mechanics / A Solids} 100 (2023) 104991.
\newblock \href {https://doi.org/10.1016/j.euromechsol.2023.104991}
  {\path{doi:10.1016/j.euromechsol.2023.104991}}.

\bibitem{KALINA_2023}
M.~Kalina, T.~Schneider, J.~Brummund, M.~Kästner, Overview of phase-field
  models for fatigue fracture in a unified framework, Engineering Fracture
  Mechanics 288 (2023) 109318.
\newblock \href {https://doi.org/10.1016/j.engfracmech.2023.109318}
  {\path{doi:10.1016/j.engfracmech.2023.109318}}.

\bibitem{KIM19961513}
J.-K. Kim, Y.-Y. Kim, Experimental study of the fatigue behavior of high
  strength concrete, Cement and Concrete Research 26~(10) (1996) 1513 -- 1523.
\newblock \href {https://doi.org/10.1016/0008-8846(96)00151-2}
  {\path{doi:10.1016/0008-8846(96)00151-2}}.

\bibitem{schneider2018untersuchungen}
S.~Schneider, J.~H{\"u}mme, S.~Marx, L.~Lohaus, {Untersuchungen zum Einfluss
  der Probek{\"o}rpergröße auf den Erm{\"u}dungswiderstand von hochfestem
  Beton}, Beton‐ und Stahlbetonbau 113~(1) (2018) 58--67.
\newblock \href {https://doi.org/10.1002/best.201700051}
  {\path{doi:10.1002/best.201700051}}.

\bibitem{Baktheer_2022_tip_bearing}
A.~Baktheer, H.~Spartali, R.~Chudoba, J.~Hegger, Concrete splitting and
  tip-bearing effect in the bond of anchored bars tested under fatigue loading
  in the push-in mode: An experimental investigation, Materials and Structures
  55~(3) (2022) 101.
\newblock \href {https://doi.org/10.1617/s11527-022-01935-7}
  {\path{doi:10.1617/s11527-022-01935-7}}.

\bibitem{BECKS_2024}
H.~Becks, M.~Classen, New insights into the load sequence effect: Experimental
  characterization and incremental modeling of plain high-strength concrete
  under mode {II} fatigue loading with variable amplitude, International
  Journal of Fatigue 185 (2024) 108334.
\newblock \href {https://doi.org/10.1016/j.ijfatigue.2024.108334}
  {\path{doi:10.1016/j.ijfatigue.2024.108334}}.

\bibitem{grzybowski1993damage}
M.~Grzybowski, C.~Meyer, Damage accumulation in concrete with and without fiber
  reinforcement, Materials Journal 90~(6) (1993) 594--604.
\newblock \href {https://doi.org/10.14359/4438} {\path{doi:10.14359/4438}}.

\bibitem{shah1984predictions}
S.~Shah, Predictions of cumulative damage for concrete and reinforced concrete,
  Mat{\'e}riaux et Construction 17~(1) (1984) 65--68.
\newblock \href {https://doi.org/10.1007/BF02474059}
  {\path{doi:10.1007/BF02474059}}.

\bibitem{jinawath1974cumulative}
P.~Jinawath, Cumulative fatigue damage of plain concrete in compression, Ph.D.
  thesis, University of Leeds (1974).

\bibitem{becks_mode_II}
H.~Becks, M.~Classen, Mode ii behavior of high-strength concrete under
  monotonic, cyclic and fatigue loading, Materials 14~(24) (2021).
\newblock \href {https://doi.org/10.3390/ma14247675}
  {\path{doi:10.3390/ma14247675}}.

\bibitem{Becks_2022_Characterization}
H.~Becks, M.~Aguilar, R.~Chudoba, M.~Classen, Characterization of high-strength
  concrete under monotonic and fatigue mode {II} loading with actively
  controlled level of lateral compression, Materials and Structures 55~(10)
  (2022) 252.
\newblock \href {https://doi.org/10.1617/s11527-022-02087-4}
  {\path{doi:10.1617/s11527-022-02087-4}}.

\bibitem{Aguilar_framcos_2023}
M.~Aguilar, A.~Baktheer, H.~Becks, M.~Classen, R.~Chudoba, Fatigue-induced
  concrete fracture under combined compression and shear studied using standard
  cylinder and refined punch-through shear test setup, 11th International
  Conference on Fracture Mechanics of Concrete and Concrete Structures (2023).
\newblock \href {https://doi.org/10.21012/FC11.092350}
  {\path{doi:10.21012/FC11.092350}}.

\bibitem{AGUILAR_2024}
M.~Aguilar, A.~Baktheer, R.~Chudoba, Multi-axial fatigue of high-strength
  concrete: Model-enabled interpretation of punch-through shear test response,
  Engineering Fracture Mechanics 311 (2024) 110532.
\newblock \href {https://doi.org/10.1016/j.engfracmech.2024.110532}
  {\path{doi:10.1016/j.engfracmech.2024.110532}}.

\end{thebibliography}

\end{document}